\documentclass[12pt]{article}

\title{Galois Theory for Braided Tensor Categories \\ and the Modular Closure}
\author{Michael M\"uger\thanks{Supported by EU TMR Network `Noncommutative Geometry'.} \\ Dipartimento di Matematica, Universit\`{a} di Roma ``Tor Vergata''\\ Via della Ricerca Scientifica, 00133 Roma, Italy\\ Email: mueger@axp.mat.uniroma2.it}
\date{January 28, 1999}

\newlength{\dinwidth}
\newlength{\dinmargin}
\usepackage{amsgen,amsfonts,amssymb,tan-gle}

\def\mobj#1{\raise .4\unitlens\hbox{\put(0,0){$#1$}}}

\def\1#1{{\bf #1}}
\def\2#1{{\cal #1}}
\def\3#1{{\sl #1}}
\def\4#1{{\tt #1}}
\def\5#1{{\sf #1}}
\def\6#1{{\mathfrak #1}}
\def\7#1{{\mathbb #1}}

\newcommand{\be}{\begin{equation}}
\newcommand{\ee}{\end{equation}}
\newcommand{\ba}{\begin{array}}
\newcommand{\ea}{\end{array}}
\newcommand{\bea}{\begin{eqnarray}}
\newcommand{\eea}{\end{eqnarray}}
\newcommand{\bean}{\begin{eqnarray*}}
\newcommand{\eean}{\end{eqnarray*}}
\newcommand{\nn}{\nonumber}

\newcommand{\ve}{\varepsilon}

\newcommand{\impl}{\Rightarrow}

\newcommand{\restr}{\upharpoonright}
\newcommand{\ol}{\overline}
\newcommand{\ul}{\underline}

\newcommand{\id}{\mbox{id}}
\newcommand{\obj}{\mbox{Obj}}
\newcommand{\mcirc}{\,\circ\,}
\newcommand{\gal}{\mbox{Gal}}
\newcommand{\aut}{\mbox{Aut}}

\newcommand{\qed}{\hfill$\blacksquare$}

\newcommand{\qfts}{quantum field theories}
\newcommand{\poinc}{Poincar\'{e}}

\newcommand{\npb}{Nucl. Phys. \1B}
\newcommand{\cmp}{Commun. Math. Phys. }

\newtheorem{defin}{Definition}[section]
\newtheorem{lemma}[defin]{Lemma}
\newtheorem{prop}[defin]{Proposition}
\newtheorem{theorem}[defin]{Theorem}
\newtheorem{coro}[defin]{Corollary}
\newtheorem{conj}[defin]{Conjecture}

\newcommand{\bdefin}{\begin{defin}}
\newcommand{\blemma}{\begin{lemma}}
\newcommand{\bprop}{\begin{prop}}
\newcommand{\btheor}{\begin{theorem}}
\newcommand{\bcoro}{\begin{coro}}
\newcommand{\edefin}{\end{defin}}
\newcommand{\elemma}{\end{lemma}}
\newcommand{\eprop}{\end{prop}}
\newcommand{\etheor}{\end{theorem}}
\newcommand{\ecoro}{\end{coro}}
\newcommand{\bconj}{\begin{conj}}
\newcommand{\econj}{\end{conj}}
\newcommand{\Hom}{\mbox{Hom}}

\newcommand{\prf}{{\it Proof. }}
\newcommand{\rem}{{\it Remark. }}
\newcommand{\rems}{{\it Remarks. }}
\newcommand{\sectreset}[1]{\section{#1}\setcounter{equation}{0}}

\begin{document}
\maketitle\noindent

\abstract{Given a braided tensor $*$-category $\2C$ with conjugate (dual) objects and
irreducible unit together with a full symmetric subcategory $\2S$ we define a crossed
product $\2C\rtimes\2S$. This construction yields a tensor $*$-category with conjugates
and an irreducible unit. (A $*$-category is a category enriched over $\mbox{Vect}_\7C$ 
with positive $*$-operation.) A Galois correspondence is established between intermediate
categories sitting between $\2C$ and $\2C\rtimes\2S$ and closed subgroups of the Galois 
group $\gal(\2C\rtimes\2S/\2C)=\aut_\2C(\2C\rtimes\2S)$ of $\2C$, the latter being 
isomorphic to the
compact group associated to $\2S$ by the duality theorem of Doplicher and Roberts. 
Denoting by $\2D\subset\2C$ the full subcategory of degenerate objects, i.e.\ objects 
which have trivial monodromy with all objects of $\2C$, the braiding of $\2C$ extends to
a braiding of $\2C\rtimes\2S$ iff $\2S\subset\2D$. Under this condition $\2C\rtimes\2S$ 
has no non-trivial degenerate objects iff $\2S=\2D$. If the original category $\2C$ is 
rational 
(i.e.\ has only finitely many isomorphism classes of irreducible objects) then the same 
holds for the new one. The category $\ol{\ol{\2C}}\equiv\2C\rtimes\2D$ is called the 
{\it modular closure} of $\2C$ since in the rational case it is modular, i.e.\ gives 
rise to a unitary representation of the modular group $SL(2,\7Z)$. (In passing we prove 
that every braided tensor $*$-category with conjugates automatically is a ribbon 
category, i.e.\ has a twist.) If all simple objects of $\2S$ have dimension one the
structure of the category $\2C\rtimes\2S$ can be clarified quite explicitly in terms of 
group cohomology.}

\sectreset{Introduction}
Since in this paper we are concerned with symmetric and braided tensor (or monoidal)
categories \cite{cwm} it may be useful to sketch the origin of some of the pertinent 
ideas. Symmetric
tensor categories were formalized in the early sixties, but they are implicit in the 
earlier Tannaka-Krein duality theory for compact groups. Further motivation for their 
analysis came from Grothendieck's theory of motives and led to Saaveda Rivano's work 
\cite{sr}, which was corrected and extended in \cite{del1}. These formalisms reconstruct
a group (compact topological in the first, algebraic in the second case) from the
category of its representations, the latter being {\it concrete}, i.e.\ consisting of 
vector (Hilbert) spaces and linear maps between these. 

In the operator algebraic approach to quantum field theory it was realized around 1970 
that the category of localized superselection sectors ($\cong$ physically relevant 
representations of the $C^*$-algebra $\2A$ of observables) is symmetric monoidal, 
cf.\ \cite{r2}. This category being a category
of endomorphisms of $\2A$ -- not of vector spaces -- the existing duality theorems did
not apply. This led Doplicher and Roberts to developing their characterization \cite{dr1}
of representation categories of compact groups as {\it abstract} symmetric tensor 
categories satisfying certain additional axioms. This 
result allowed the solution \cite{dr2} of the longstanding problem of (re-)constructing
a net of charged field algebras $\2F$ which intertwine the inequivalent representations 
of $\2A$ and have nice properties like Bose-Fermi commutation relations. (In fact,
{\it assuming} the duality theorem for abstract symmetric categories such a 
reconstruction result existed much earlier \cite{r1}.) At the same time and independently
Deligne extended \cite{del2} the earlier works \cite{sr,del1} by identifying a necessary
and sufficient condition for an abstract symmetric tensor category to be the 
representation category of an algebraic group. The crucial property is that all objects 
in the given category have integer dimension. (For symmetric $C^*$-tensor categories 
this is automatic \cite[Cor. 2.15]{dr1} as a consequence of Hilbert space positivity.)

Braided tensor categories without the symmetry requirement entered the scene only in
the eighties. From a theoretical point of view braided tensor categories
are most naturally `explained' by identifying \cite{js} them as 2-categories with
tensor product and only one object, which in turn are just 3-categories with only 
one object and one 1-morphism. (All these notions are easiest to deal with in the
strict case, which for (symmetric, braided) tensor categories does not imply a loss of
generality in view of the coherence theorems \cite{cwm,js}.)
But the main reason for their recent prominence is their relation to certain algebraic
structures arising in physics (Yang-Baxter equation, quantum groups) and to topological 
invariants of knots, links and 3-manifolds. The latter subject 
was boosted by V.\ Jones' construction of a new knot invariant which was soon discovered
to be related to the quantum group $SU_q(2)$ where $q$ is a root of unity, and
subsequently invariants of 3-manifolds were constructed for all quantum groups
at roots of unity. The theory reached a certain state of maturity when it was understood
that the crucial ingredient underlying these invariants of 3-manifolds is a certain class
of braided tensor categories which are called modular \cite{t}. A modular category is
a braided tensor category which (i) has a twist \cite{t} or balancing \cite{js}, (ii) is
rational -- i.e.\ has only finitely many isomorphism classes of irreducible objects --
 and (iii) non-degenerate. Here non-degeneracy means that an irreducible object
$\rho$ for which $\ve(\rho,\sigma)\circ\ve(\sigma,\rho)=\id_{\sigma\rho}\ \forall\sigma$
is equivalent to the unit object $\iota$. (The designation of such categories as modular
is owed to the fact that they give rise to a finite dimensional representation of the 
modular group $SL(2,\7Z)$ \cite{t}, see also \cite{khr1}.)
The role of the quantum groups then reduces just to providing several infinite families
of modular categories (roughly, one for every pair (root of unity, classical Lie 
algebra)). Another construction of modular categories starts from link invariants, cf.\
\cite[Chap.\ XII]{t}, \cite{tw}. Finally, braided tensor categories appear naturally 
also in the superselection theory of quantum field theories in low dimensional 
spacetimes, cf.\ eg.\ \cite{kmr}. In many cases, as for the WZW and orbifold models, 
these categories actually turn out to be modular. Let $\2A$ be a quantum field theory in
$1+1$ dimensions and let $\2C$ be the braided category of superselection sectors with 
finite statistics. Since the full subcategory $\2D\subset\2C$ of degenerate sectors is 
symmetric, the Doplicher-Roberts construction can be applied to $\2A$ and $\2D$ and 
yields new theory $\2F$. In \cite{khr1} Rehren conjectured that the representation 
category of $\2F$ is non-degenerate. Under the assumption that $\2A$ has only finitely 
many irreducible degenerate sectors
this was proved by the author in \cite{mue4}. The aim of the present paper is to give a 
purely categorical analogue of this construction (without the finiteness restriction).

More precisely, given a braided tensor category $\2C$ which is enriched over 
${\bf \mbox{Vect}_\7C}$,
has a positive $*$-operation, conjugate (dual) objects, direct sums and subobjects and 
an irreducible unit object together with a symmetric subcategory $\2S$ satisfying the 
same properties, we define a crossed product $\2C\rtimes\2S$. (The existence of direct 
sums and subobjects (in the sense of \cite{glr}) is no serious restriction since it can 
always be achieved by embedding the category in a bigger one \cite[Appendix]{lr}.) 
This construction proceeds in two steps. First we define a  tensor category 
$\2C\rtimes_0\2S$ which has the same objects and tensor product as $\2C$ but bigger 
spaces of arrows, i.e.\ 
\be \Hom_{\2C\rtimes_0\2S}(\rho,\sigma)\supset \Hom_{\2C}(\rho,\sigma)\quad\forall 
  \rho,\sigma\in\2C. \ee
Of course, we have to prove that $\2C\rtimes_0\2S$ satisfies all axioms of a tensor
$*$-category. The new category inherits the braiding $\ve$ from $\2C$ iff $\2S$ contains 
only degenerate objects, thus iff $\2S\subset\2D$ where $\2D\subset\2C$ is the full 
subcategory of degenerate objects as above.
(If this condition is not fulfilled naturality of $\ve$ fails for 
some of the new morphisms of $\2C\rtimes_0\2S$). Now, $\2C\rtimes_0\2S$ will be closed 
under direct sums, but usually not under subobjects. Thus we apply the above-mentioned
procedure of \cite{lr} in order to obtain a category $\2C\rtimes\2S$ which is closed 
under direct sums and subobjects. Then the braiding of $\2C\rtimes_0\2S$ -- if it exists
-- extends to $\2C\rtimes\2S$. The result of this construction is again a 
tensor category with positive $*$-operation and conjugates, direct sums, subobjects
and irreducible unit. $\2C\rtimes\2S$ is braided if $\2S\subset\2D$. Under this 
condition we prove that $\2C\rtimes\2S$ has no degenerate objects iff $\2S=\2D$. 
The category $\ol{\ol{\2C}}=\2C\rtimes\2D$ is called the {\it modular closure} of $\2C$ 
since in the rational case (where there are only finitely many isomorphism classes of 
irreducible objects) it is modular. (In particular, $\ol{\ol{\2C}}$ is rational if $\2C$
is.) The modular closure $\ol{\ol{\2C}}$ is non-trivial, i.e.\ has irreducible objects 
which are not equivalent to the unit, iff $\2C$ is not symmetric, thus not completely 
degenerate. Define the absolute Galois group $\gal(\2C)$ of a braided tensor category 
$\2C$ to be the compact group associated to the symmetric tensor category $\2D(\2C)$ by 
the duality theorem of Doplicher and Roberts. 
For every symmetric category $\2S\subset\2C$ we establish a Galois correspondence 
between subcategories $\2B$ of $\2C\rtimes\2S$ containing $\2C$ and closed 
subgroups $H$ of the relative Galois group $G=\aut_\2C(\2C\rtimes\2S)\cong\gal(\2S)$, 
given by $\2B=(\2C\rtimes\2S)^H$ and $H=\aut_\2B(\2C\rtimes\2S)$. The normal subgroups 
$H$ correspond to the subextensions $\2C\rtimes\2T$ where $\2T\subset\2S$ and 
$\gal(\2T)\cong G/H$. If $\2S\subset\2D$ then $\2C\rtimes\2S$ is a braided subextension
of $\ol{\ol{\2C}}=\2C\rtimes\2D$, the absolute Galois group $\gal(\2C\rtimes\2S)$
being isomorphic to $H=\aut_{\2C\rtimes\2S}(\ol{\ol{\2C}})$.
Giving an explicit description of the (isomorphism classes of) irreducible objects of
$\2C\rtimes\2S$ is difficult in general, but if all irreducible objects of $\2S$ have 
dimension one, corresponding to abelian $\gal(\2S)$, the structure of the category 
$\2C\rtimes\2S$ can be clarified quite explicitly in terms of group cohomology.

We briefly describe the organization of the paper.
In Sect.\ 2 we give precise definitions and several preparatory results on $C^*$-tensor
categories. In particular we prove that they are automatically ribbon categories, i.e.\
have a twist. In Sect.\ 3 the crossed product $\2C\rtimes\2S$ is defined and proved to
be a $C^*$-tensor category. Then, in Sect.\ 4 we prove that $\2C\rtimes\2D$ is 
non-degenerate and establish the Galois correspondence. In Sect.\ 5 we enlarge on 
abelian extensions, the case of supergroups and make some further remarks on the case 
$\2S\not\subset\2D$.

\sectreset{Definitions and Preparations}
\subsection{Some Results on $C^*$-Tensor Categories}
We begin by establishing our notation concerning tensor categories.
Objects will be denoted by small Greek letters $\rho,\sigma$, etc. The set of arrows 
(morphisms) between $\rho$ and $\sigma$ in the category $\2C$ is $\Hom_\2C(\rho,\sigma)$,
where the subscript $\2C$ is omitted when there is no danger of confusion. The
identity arrow of $\rho$ is $\id_\rho$, and composition of arrows is denoted by $\circ$.
The tensor product $\rho\otimes\sigma$ of objects will abbreviated by $\rho\sigma$. All 
tensor categories in this paper are supposed small and strict, thus when we mention
these conditions it is only for emphasis. (A tensor category is strict if the tensor
product satisfies associativity $\rho(\sigma\eta)=(\rho\sigma)\eta$ `on the nose' and
there is a unit object $\iota$ satisfying $\rho\iota=\iota\rho=\rho\ \forall\rho$.)
Given two arrows $R\in\Hom(\rho,\sigma), R'\in\Hom(\rho',\sigma')$ there is an arrow 
$R\times R'\in\Hom(\rho\rho',\sigma\sigma')$. The mapping $(R, R')\mapsto R\times R'$
is associative, satisfies $\id_\iota\times R=R\times\id_\iota=R$ and the interchange law
\be (S\circ R)\times(S'\circ R')=S\times S'\circ R\times R' \ee
if $S\in\Hom(\sigma,\tau), S'\in\Hom(\sigma',\tau')$. A tensor category $\2C$ is braided 
if there is a family of invertible arrows 
$\{ \ve(\rho,\sigma)\in\Hom(\rho\sigma,\sigma\rho),\ \rho,\sigma\in\2C \}$, natural in 
both variables and satisfying
\bea \ve(\rho,\sigma_1\sigma_2) &=& \id_{\sigma_1}\times\ve
   (\rho,\sigma_2)\mcirc \ve(\rho,\sigma_1)\times\id_{\sigma_2}, \\
   \ve(\rho_1\rho_2,\sigma)&=&\ve(\rho_1,\sigma)\times\id_{\rho_2}\mcirc
   \id_{\rho_1}\times\ve(\rho_2,\sigma) \eea
for all $\rho_i,\sigma_i$. A braided tensor category is symmetric if the braiding
satisfies $\ve(\rho,\sigma)\circ\ve(\sigma,\rho)=\id_{\sigma\rho}\ \forall\rho,\sigma$.

All categories in this paper will be enriched over ${\bf \mbox{Vect}_\7C}$, but we do not
presuppose familiarity with this notion. A {\it complex tensor category} is a tensor 
category, for which the sets $\Hom(\rho,\sigma)$ of arrows are complex vector spaces and the 
composition $\circ$ and tensor product $\times$ of arrows are bilinear. 
A $*$-operation on a complex tensor category is map which assigns to an arrow
$X\in\Hom(\rho, \sigma)$ another arrow $X^*\in\Hom(\sigma,\rho)$. This map has to be antilinear, 
involutive ($X^{**}=X$), contravariant ($(S\circ T)^*=T^*\circ S^*$) and 
monoidal ($(S\times T)^*=S^*\times T^*$). A $*$-operation is {\it positive} iff 
$X^*\circ X=0$ implies $X=0$. A {\it tensor $*$-category} is a complex tensor category
with a positive $*$-operation. For such categories we admit only unitary braidings.

An object $\rho$ is called irreducible or simple if $\Hom(\rho,\rho)=\7C\,\id_\rho$.
As usual, two objects $\rho,\sigma$ are equivalent (or isomorphic) iff $\Hom(\sigma,\rho)$ 
contains an invertible arrow. In $W^*$-categories $\Hom(\sigma,\rho)$ then contains a 
unitary by polar decomposition of morphisms \cite[Cor.\ 2.7]{glr}.
An object $\sigma$ is a {\it subobject} of 
$\rho$, denoted $\sigma\prec\rho$, iff $\Hom(\sigma,\rho)$ contains an isometry. 
Note that this notion of subobjects differs from the standard one of category theory 
\cite{cwm}, cf.\ also the remarks in \cite[p.\ 98]{glr}. A tensor $*$-category is closed 
under subobjects (or: has subobjects) if for every orthogonal projection 
$E\in\Hom(\rho,\rho)$ there is an object $\sigma$ and an isometry $V\in\Hom(\sigma,\rho)$
such that $V\circ V^*=E$. A tensor $*$-category has (finite) {\it direct sums} iff for 
every pair 
$\rho_1,\rho_2$ there are $\tau$ and isometries $V_i\in\Hom(\rho_i,\tau)$ such that 
$V_1\circ V_1^*+V_2\circ V_2^*=\id_\tau$. Then we write $\tau\cong\rho_1\oplus\rho_2$. 
Note that every $\tau'\cong\tau$ is a direct sum of $\rho_1,\rho_2$, too.
A tensor $*$-category can always canonically be extended to a tensor $*$-category
with direct sums and subobjects \cite[Appendix]{lr}.

From now on all categories are tensor $*$-categories. In the present setting it is 
convenient to define conjugate (dual) objects in a way which differs
slightly from the one for rigid (autonomous) tensor categories \cite{js,t}.
We give only the main definitions and facts and refer to \cite{lr} for the details.
An object $\ol{\rho}$ is said to be conjugate to $\rho$ if there 
are $R\in\Hom(\iota,\ol{\rho}\rho), \ol{R}\in\Hom(\iota,\rho\ol{\rho})$ satisfying the 
{\it conjugate equations}:
\be \ol{R}^*\times\id_\rho\mcirc\id_\rho\times R=\id_\rho, 
  \quad\quad  R^*\times\id_{\ol{\rho}}\mcirc\id_{\ol{\rho}}\times
  \ol{R}= \id_{\ol{\rho}}. \label{conj-eq0}\ee
A category $\2C$ has conjugates if every object $\rho\in\2C$ has a conjugate 
$\ol{\rho}\in\2C$. If $\rho$ is irreducible, then an irreducible conjugate
$\ol{\rho}$ is unique up to isomorphism and (upon proper normalization of $R, \ol{R}$)
$R^*\circ R=\ol{R}^*\circ\ol{R}\in\Hom(\iota,\iota)$ is independent of the choice of
$R, \ol{R}$. Then the dimension defined via $d(\rho)\id_\iota=R^*\circ R$ is in
$[1,\infty)$ and satisfies $d(\rho)=d(\ol{\rho})$. For reducible $\rho$ we 
admit only {\it standard solutions} \cite{lr} of (\ref{conj-eq0}). This means that
$R_\rho=\sum_i \ol{W}_i\times W_i\mcirc R_i$ where $\rho=\oplus_i \rho_i$ is a
decomposition into irreducibles effected by the isometries $W_i$ and $R_i$ is (part of)
a normalized solution of (\ref{conj-eq0}) for $\rho_i$. Then the definition 
$d(\rho)=R_\rho^*\circ R_\rho$ extends to reducible objects and yields a multiplicative
dimension function. (This dimension is subject to the same restriction as the Jones 
index of an inclusion of factors, cf.\ \cite{lr}. Note that the braiding does not play a
role here, yet the categorical dimension coincides with the $q$-dimension for
representation categories of quantum groups \cite{rt}.)

The more specific notion of $C^*$-tensor categories will not be needed explicitly 
in this paper. But since we wish to make use of results of \cite{glr,dr1,lr} we will
prove that many tensor $*$-categories are automatically $C^*$-tensor categories.
Now, a $C^*$-tensor category is a complex tensor category with a 
$*$-operation. Furthermore, the spaces $\Hom(\rho,\sigma),\ \rho,\sigma\in\2C$ are
Banach spaces and the norms satisfy
\bea \| Y\circ X\| & \le &\|X\|\,\|Y\|, \label{subm} \\
   \| X^*\circ X\| & = &\|X\|^2 \label{cstar0} \eea
for $X\in\Hom(\rho,\sigma), Y\in\Hom(\sigma,\eta)$. (Then the algebras 
$\Hom(\rho,\rho),\,\rho\in\2C$ are $C^*$-algebras.) See the cited references for
examples.

It is known \cite{lr} that in a $C^*$-tensor category with conjugates and an irreducible
unit, i.e.\ $\Hom(\iota,\iota)=\7C\,\id_\iota$, all spaces of arrows are finite 
dimensional. The following result is a converse, which generalizes a well-known fact on 
finite dimensional $\7C$-algebras.
\bprop Let $\2C$ be a $\mbox{Vect}_\7C^{fin}$-category, i.e.\ a category 
where $\Hom(\rho,\sigma)$ is a finite dimensional $\7C$-vector space for every pair 
$\rho,\sigma\in\2C$, the composition $\circ$ being bilinear.
Then $\2C$ is a $C^*$-category iff there is a positive $*$-operation. \label{cstar}\eprop
\prf If $\2C$ is a $C^*$-category there is a $*$-operation by definition. 
Positivity follows from (\ref{cstar0}). Assume conversely the existence of a positive
$*$-operation. In particular, $*$ gives rise to a positive involution on the algebras 
$\Hom(\rho, \rho),\ \rho\in\2C$. The latter being 
finite dimensional $\7C$-algebras, this implies semisimplicity and the existence of 
unique $C^*$-norms. Now we consider the $*$-algebras $M(\rho_1,\ldots,\rho_n)$
\cite[p.\ 86]{glr} associated to $n$ objects (which, roughly speaking, are the algebras
generated by the arrows between the objects $\rho_1,\ldots,\rho_n$). For an element 
$\hat{X}=(X_{ij})$ of $M(\rho_1,\ldots,\rho_n),\ \ \hat{X}^*\hat{X}=0$ is 
equivalent to $X_{ij}^*\circ X_{ij}=0\ \forall i,j=1,\ldots,n$.
Since by assumption this holds only if all $X_{ij}$ vanish, the $*$-involution 
of $M(\rho_1,\ldots,\rho_n)$ is positive and also $M(\rho_1,\ldots,\rho_n)$ is a 
$C^*$-algebra. Now we define the norm on $\Hom(\rho, \sigma)$ by 
\be \|X\| = \sqrt{\|X^*\circ X\|}, \quad X\in\Hom(\rho, \sigma),\label{norm}\ee
where the norm on the right hand side is the one of $M(\rho,\sigma)$.
Since the algebras $M$ form a directed system the norm of $X^*\circ X$ is the same in, 
say, $M(\rho,\sigma,\eta)$ and thus well defined. As an immediate consequence we have
$\|X\| = \|X^*\|$, and the submultiplicativity of the norms
\be \| Y\circ X\| \le \|X\|\,\|Y\| \ee
for $X\in\Hom(\rho,\sigma), Y\in\Hom(\sigma,\eta)$ required of a $C^*$-category follows
from submultiplicativity in $M(\rho,\sigma,\eta)$. The $C^*$-condition (\ref{cstar0})
follows from the $C^*$-property of $M(\rho,\sigma)$.  \qed\\
\rem This result is probably well known among experts, but to the best of the author's
knowledge it never appeared in print. Yet it is used implicitly in \cite{w} where
certain categories are proved to have a positive $*$-operation and concluded to be
$C^*$-categories.

In the above result we did not assume irreducibility of the unit $\iota$, viz.
$\Hom(\iota,\iota)=\7C\,\id_\iota$. From now on all categories in this paper will be 
assumed to have this property, which has been called connectedness \cite{baez}. We will 
remark on the disconnected case in the outlook.

We summarize the properties of the categories we will study.
\bdefin A T$C^*$ is a small strict tensor $*$-category with conjugates, direct sums, 
subobjects, finite dimensional spaces of arrows and an irreducible unit object. 
A BT$C^*$ is a T$C^*$ with a unitary braiding. A ST$C^*$ is a symmetric BT$C^*$. \edefin
\rem All concepts in this definition which are not standard category theory are from 
\cite{glr,dr1,lr}. That they were arrived at independently under the name `unitary
categories' \cite[Sect.\ II.5]{t} underlines their naturality.

In the literature on braided tensor categories additional pieces of structure have been 
considered, mostly motivated by the study of topological invariants of 3-manifolds.
\bdefin A twist \cite{t} or balancing \cite{js} for a braided tensor category $\2C$ is 
a family $\{ \kappa(\rho)\in\Hom(\rho,\rho),\ \rho\in\2C\}$ of invertible arrows
satisfying naturality
\be T\circ\kappa(\rho)=\kappa(\sigma)\circ T \quad\forall T\in\Hom(\rho,\sigma)
\label{tw1}\ee
and the conditions 
\be \kappa(\rho_1\rho_2)=\kappa(\rho_1)\times\kappa(\rho_1)\mcirc\ve(\rho_2,\rho_1)
   \mcirc\ve(\rho_1,\rho_2) \quad \forall \rho_1,\rho_2, \label{tw2}\ee
\be \kappa(\ol{\rho})\times\id_\rho\mcirc R=\id_{\ol{\rho}}\times\kappa(\rho)\mcirc R
\label{tw3}\ee
for every standard solution $(\rho,\ol{\rho}, R, \ol{R})$ of the conjugate equations. 
In a tensor $*$-category $\kappa(\rho)$ is required to be unitary. \edefin
\rems 1. The condition (\ref{tw1}) is equivalent to saying the $\kappa$ is a natural 
transformation of the identity functor to itself. (The set of these was called 
the center of $\2C$ in \cite{glr}.) \\
2. In \cite{js} the definition of a twist does not include (\ref{tw3}). There a 
category with conjugates and a twist satisfying (\ref{tw3}) is called tortile.\\
3. If $\rho$ is irreducible then we define $\omega(\rho)\in\7C$ via
$\kappa(\rho)=\omega(\rho)\id_\rho$. 

\begin{figure}[t]
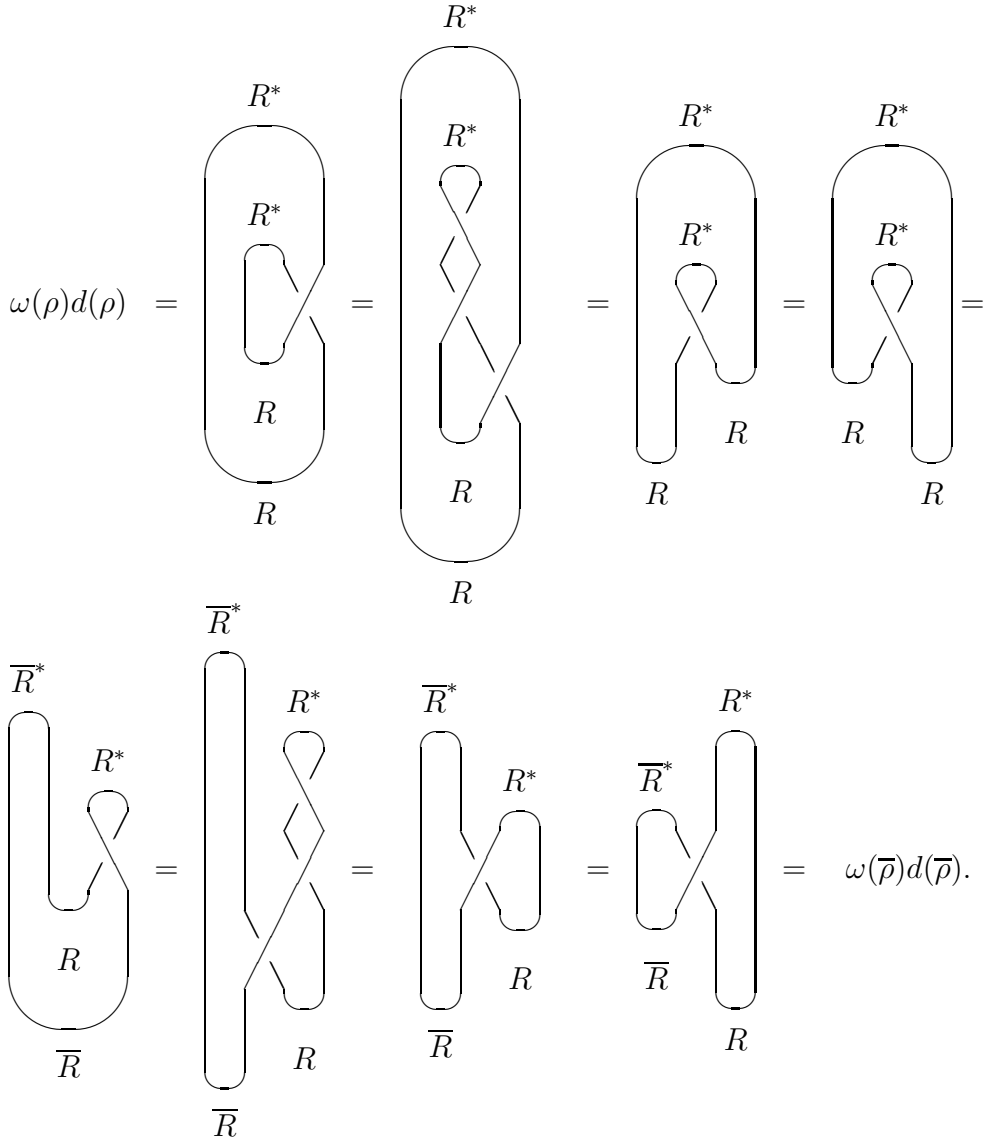

\[\ba{ccccccccc}
\omega(\rho)d(\rho) &=& 
\begin{tangle}
\step[1.5]\object{R^*}\\
\mcoev\\
\hh\id\step[1.5]\object{R^*}\step[1.5]\id\\
\hh\id\step\hcoev\step\id\\
\id\step\id\step\hxx\\
\hh\id\step\hev\step\id\\
\hh\id\step[1.5]\object{R}\step[1.5]\id\\
\mev\\
\step[1.5]\object{R}
\end{tangle}
& = &
\begin{tangle}
\step[1.5]\object{R^*}\\
\mcoev\\
\hh\id\step[1.5]\object{R^*}\step[1.5]\id\\
\hh\id\step\hcoev\step\id\\
\id\step\hx\step\id\step\\
\id\step\hxx\step\id\step\\
\id\step\id\step\hxx\\
\hh\id\step\hev\step\id\\
\hh\id\step[1.5]\object{R}\step[1.5]\id\\
\mev\\
\step[1.5]\object{R}
\end{tangle}
& = &
\begin{tangle}
\step[1.5]\object{R^*}\\
\mcoev\\
\hh\id\step[1.5]\object{R^*}\step[1.5]\id\\
\hh\id\step\hcoev\step\id\\
\id\step\hx\step\id\\
\hh\id\step\id\step\hev\\
\hh\id\step\id\step[1.5]\object{R}\\
\hh\hev\\
\hstep\object{R}
\end{tangle}
& = &
\begin{tangle}
\step[1.5]\object{R^*}\\
\mcoev\\
\hh\id\step[1.5]\object{R^*}\step[1.5]\id\\
\hh\id\step\hcoev\step\id\\
\id\step\hx\step\id\\
\hh\hev\step\id\step\id\\
\hh\hstep\object{R}\step[1.5]\id\step\id\\
\hh\Step\hev\\
\step[2.5]\object{R}
\end{tangle}=
\\
\begin{tangle}
\hstep\object{\ol{R}^*}\\
\hh\hcoev\\
\hh\id\step\id\step[1.5]\object{R^*}\\
\hh\id\step\id\step\hcoev\\
\id\step\id\step\hx\\
\hh\id\step\hev\step\id\\
\hh\id\step[1.5]\object{R}\step[1.5]\id\\
\mev\\
\step[1.5]\object{\ol{R}}
\end{tangle} 
& = &
\begin{tangle}
\hstep\object{\ol{R}^*}\\
\hh\hcoev\\
\hh\id\step\id\step[1.5]\object{R^*}\\
\hh\id\step\id\step\hcoev\\
\id\step\id\step\hx\\
\id\step\id\step\hxx\\
\id\step\hxx\step\id\\
\hh\id\step\id\step\hev\\
\hh\id\step\id\step[1.5]\object{R}\\
\hh\hev\\
\hstep\object{\ol{R}}
\end{tangle} 
& = &
\begin{tangle}
\hstep\object{\ol{R}^*}\\
\hh\hcoev\\
\hh\id\step\id\step[1.5]\object{R^*}\\
\hh\id\step\id\step\hcoev\\
\id\step\hxx\step\id\\
\hh\id\step\id\step\hev\\
\hh\id\step\id\step[1.5]\object{R}\\
\hh\hev\\
\hstep\object{\ol{R}}
\end{tangle} 
& = &
\begin{tangle}
\step[2.5]\object{R^*}\\
\hh\Step\hcoev\\
\hh\hstep\object{\ol{R}^*}\step[1.5]\id\step\id\\
\hh\hcoev\step\id\step\id\\
\id\step\hxx\step\id\\
\hh\hev\step\id\step\id\\
\hh\hstep\object{\ol{R}}\step[1.5]\id\step\id\\
\hh\Step\hev\\
\step[2.5]\object{R}
\end{tangle} 
& = &
\omega(\ol{\rho})d(\ol{\rho}).
\ea\]
\caption{Proof of $\omega(\rho)=\omega(\ol{\rho})$}
\label{kappa}\end{figure}

A remarkable feature of BT$C^*$s is that they automatically possess a canonically defined
twist. It is defined and studied in \cite[Thm.\ 4.2]{lr}, where however the property
(\ref{tw3}) was not proved. 
\bprop BT$C^*$s are ribbon categories, i.e.\ have a twist. \label{twist}\eprop
\prf In \cite[Sect.\ 4]{lr} for every BT$C^*\ \2C$ a family 
$\{ \kappa(\rho)\in\Hom(\rho,\rho),\ \rho\in\2C\}$ satisfying (\ref{tw1},\ref{tw2}) was
defined, the $\kappa$'s being unitary whenever the braiding $\ve$ is unitary. (Recall 
that we assume this throughout.) Thus it only remains to prove (\ref{tw3}) and in view 
of the naturality of the twist it is sufficient to consider only irreducible $\rho$, 
where (\ref{tw3}) reduces to $\omega(\rho)=\omega(\ol{\rho})$. This is done in 
Fig.\ \ref{kappa}. In the first and last equalities we have used that for $\rho$ 
irreducible and $\Hom(\rho,\rho)\ni T=C\,\id_\rho$ we have 
\be R^*\mcirc \id_{\ol{\rho}}\times T\mcirc R=\ol{R}^*\mcirc T\times\id_\rho\mcirc\ol{R}
  = C\,d(\rho). \label{bla}\ee
(Here we use $d(\rho)=R^*\circ R=\ol{R}^*\circ\ol{R}=d(\ol{\rho})$, cf.\ \cite{lr}.)
That the two ways of closing the loop in (\ref{bla}) yield the same result is used in the
fifth equality of the above calculation. The other steps use nothing more than the 
interchange law. \qed\\
\rem This argument has been adapted from algebraic quantum field theory, cf.\ 
\cite[Lemma II.5.14]{kmr}, but see also \cite{baez}. L.\ Tuset independently arrived at
essentially the same proof.

\subsection{The Galois Group of a Braided Tensor Category}
\bdefin The monodromy of two objects of a braided tensor category $\2C$ is 
\be \ve_M(\rho,\sigma)\equiv\ve(\sigma,\rho)\circ\ve(\rho,\sigma)
   \in\Hom(\rho\sigma,\rho\sigma). \ee
An object $\rho\in\2C$ is {\it degenerate} iff 
\be \ve_M(\rho,\eta)=\id_{\rho\eta}\quad\forall\eta\in\2C. \label{degen}\ee
A braided tensor category is degenerate if there is an irreducible degenerate object
which is not isomorphic to the unit object $\iota$. \edefin
\rem Clearly, a braided tensor category is symmetric iff all objects are degenerate.

\bdefin Let $\2C$ be a BT$C^*$. Then $\2D(\2C)$ is the full subcategory whose objects 
are the degenerate objects of $\2C$. \edefin

\bprop $\2D(\2C)$ is a symmetric tensor category with $*$-operation, conjugates, direct
sums, subobjects and finite dimensional spaces of morphisms.  \label{dc}\eprop
\prf For $\rho,\sigma,\eta\in\2C$ we have 
\bea \lefteqn{\ve_M(\rho\sigma,\eta)=\ve(\eta,\rho\sigma)\circ\ve(\rho\sigma,\eta)=} \\
   && \id_\rho\times\ve(\eta,\sigma)\mcirc\ve(\eta,\rho)\times\id_\sigma\mcirc
   \ve(\rho,\eta)\times\id_\sigma\mcirc\id_\rho\times\ve(\sigma,\eta). \nn\eea
It is easily seen that this reduces to $\id_{\rho\sigma\eta}$ if $\rho$ 
and $\sigma$ have trivial monodromy with $\eta$. Thus the set of degenerate objects is 
closed under multiplication. Now let $\rho\cong\oplus_{i\in I}\,\rho_i$, i.e.\ there are 
morphisms $V_i\in\Hom(\rho_i,\rho)$ such that $V_i^*\circ V_j=\delta_{i,j}\id_{\rho_i}$ and
$\sum_i V_i\circ V_i^*=\id_\rho$. Then by naturality of the braiding we have 
\be \ve_M(\rho,\eta)=\sum_i V_i\times\id_\eta \mcirc 
   \ve_M(\rho_i,\eta)\mcirc V_i^*\times\id_\eta, \ee
which implies that $\rho$ is degenerate iff all $\rho_i\prec\rho$ are degenerate. Thus
the set of degenerate objects is closed under direct sums and subobjects. In order to 
show that the conjugate of a degenerate object is degenerate, it is sufficient to
consider irreducible objects. The following equality is proved by the same argument as
already employed in the proof of Prop.\ \ref{twist}:
\bea\begin{tangle}
\object{\ol{\rho}}\Step\object{\ol{\eta}}\\
\xx\\
\xx\\
\object{\ol{\rho}}\Step\object{\ol{\eta}}
\end{tangle}
&=&
\begin{tangle}
\step[1.5]\object{R^*_\rho}\step[2.5]\object{\ol{\rho}}\step\object{\ol{\eta}}\\
\mcoev\step\id\step\id\\
\hh\id\step\mobj{R^*_\eta}\Step\id\step\id\step\id\\
\hh\id\step\hcoev\step\id\step\id\step\id\\
\id\step\id\step\hxx\step\id\step\id\\
\id\step\id\hstep\mobj{\eta}\hstep\hxx\step[.2]\mobj{\rho}\step[.8]\id\step\id\\
\hh\id\step\id\step\id\step\hev\step\id\\
\hh\id\step\id\step\id\step\obj{\ol{R}_\rho}\Step\id\\
\id\step\id\step\mev\\
\object{\ol{\rho}}\step\object{\ol{\eta}}\step[2.5]\object{\ol{R}_\eta}
\end{tangle}
\eea
Using this we see that $\ve_M(\rho,\sigma)=\id_{\rho\sigma}$ for all $\sigma$ implies
$\ve_M(\ol{\rho},\sigma)=\id_{\ol{\rho}\sigma}\ \forall\sigma$. $\2D(\2C)$ is a ST$C^*$,
since the braiding of $\2C$ is symmetric in restriction to the degenerate objects. \qed\\
\rem From the above it is clear that $\2D(\2C)$ is the correct object to be denoted as
the {\it center} of $\2C$. This is the analog for braided tensor categories of the
usual center of a monoid (=tensor 0-category), but it must not be confused with yet
another definition of a `center', namely the quantum double $\2Z(\2C)$ (which is a 
braided tensor category) of a tensor category $\2C$ (not necessarily braided). 

By the above result and Prop.\ \ref{cstar}, $\2D(\2C)$ satisfies the assumptions of the 
duality theorem of \cite{dr1}. We briefly summarize the principal results of \cite{dr1}.
Since every object of a symmetric tensor category $\2S$ satisfies 
$\ve(\rho,\rho)^2=\id_{\rho^2}$, the twist in a ST$C^*$ takes only the values
$\pm 1$. (In physics, objects with twist $+1$ and $-1$ are called bosons and fermions,
respectively.) For irreducible $\rho_1,\rho_2$ (\ref{tw2}) reduces to
$\kappa(\rho_1\rho_2)=\omega(\rho_1)\omega(\rho_2)\,\id_{\rho_1\rho_2}$, thus all
subobjects of $\rho_1\rho_2$ have the same twist. Therefore the objects with twist $+1$ 
generate a full subcategory $\2S_+$ which is again a BT$C^*$. We assume for a moment 
that $\2S$ is even, thus $\2S=\2S_+$. By \cite[Thm.\ 6.1]{dr1} there is a compact group
$G$ unique up to isomorphism such that $\2S\cong U(G)$ where $U(G)$ is a category
of finite dimensional unitary representations of $G$ containing representers for all
isomorphism (unitary equivalence) classes of irreducible representations of $G$. 
(Conceptually, the proof of may be considered as composed of two steps. 
First one shows that for a category with the above properties there is a symmetric 
$C^*$-tensor functor $F$, the embedding functor, into
the category $\2H$ of Hilbert spaces. $F$ is unique up to a natural transformation.
In the second step the Tannaka-Krein reconstruction theorem is applied
to the category $F(\2S)$ and shows that $F(\2S)$ is isomorphic
to a category of representations of a uniquely defined compact group $G$. But observe
that the proof in \cite{dr1} is independent of the Tannaka-Krein theory in that the
group $G$ is constructed simultaneously with the embedding.)

Since all objects in a category $U(G)$ have twist $+1$ the above result can not hold
if $\2S$ contains fermionic objects. Yet in this case the braiding in $\2S$ can be 
modified (`bosonized') such as to obtain an even BT$C^*$ $\2S'$ and a compact supergroup
$(G,k)$. Here $G$ is the compact group associated to $\2S'$ and $k$
is an element of order two in the center of $G$
such that the twist of an irreducible object in $\2S$ is the value of $k$ in the
associated representation of $G$. The group $G_+$ corresponding to $\2S_+$ is just the
quotient $G_+=G/\{e,k\}$.

\bdefin Let $\2C$ be a BT$C^*$. Then the absolute Galois group $\gal(\2C)$ is the 
compact group associated by Doplicher and Roberts to the center $\2D(\2C)$ of $\2C$. 
\label{agal}\edefin
\rem Strictly speaking, $\gal(\2C)$ is not a group but an isomorphism class of groups.
As soon as a representation functor $F:\2D(\2C)\rightarrow\2H$ has been chosen we have a
concrete group $\gal_F(\2C)$, the group of natural transformations from $F$ to itself
as first considered in \cite{sr}.

The following discussion serves only to motivate the terminology `modular closure'
of Sect.\ \ref{sect4} and may be ignored.

Given two irreducible objects $\rho,\sigma$ the number $Y_{\rho,\sigma}$ defined by
\be Y_{\rho,\sigma} \,\id_\iota = 
   R_\rho^*\times \ol{R}_\sigma^* \mcirc \id_{\ol{\rho}}\times
   \ve_M(\rho,\sigma)\times\id_{\ol{\sigma}} \mcirc R_\rho\times\ol{R}_\sigma =\ \ 
\begin{tangle}
\hcoev\Step\hcoev\\
\id\step\xx\step\id\\
\id\step[.6]\mobj{\rho}\step[.4]\xx\mobj{\sigma}\step\id\\
\hev\Step\hev
\end{tangle}
\ee
depends only on the isomorphism classes of $\rho,\sigma$. 

\bdefin A category is {\it rational} if the number of isomorphism classes of irreducible
objects is finite. \edefin
In a rational category $Y$ gives rise to a (finite) matrix indexed by the isomorphism 
classes of irreducible objects.
\bdefin A rational BT$C^*$ is modular if the matrix $Y$ is invertible. \edefin
\rem Recall that the existence of a twist which is usually required from a modular
category \cite{t} is automatic in BT$C^*$s.

\bprop A rational BT$C^*$ is modular iff it is non-degenerate. In the non-degene\-rate
case $Y$ is proportional to a unitary matrix $S$ which together with a certain matrix
$T\propto \mbox{diag}(\omega_i)$ gives rises to a unitary representation of $SL(2,\7Z)$.
\eprop
\prf The statement is the categorical version of a result from \cite{khr1} and can be 
proved by straightforward adaption of the arguments of \cite[Sect.\ 5]{khr1} to the 
framework of BT$C^*$s. (The factor $d_id_j$ in \cite[(5.11)]{khr1} is accounted for by 
the different normalizations of the $R$'s in \cite{khr1} and the present paper.)
We refrain from giving details since that would use too much space and will not be 
used in this paper. The claimed fact will be contained as a special case in a more
general result, proved in \cite{mue8}. \qed

\sectreset{Crossed Product of Braided Tensor $*$-Categories by Symmetric Subcategories}\label{sect3}
\subsection{Definition of the Crossed Product}
We assume that $\2C$ has direct sums and subobjects, which can be interpreted by saying
that reducible objects are always completely reducibility, or $\2C$ is semi-simple.
This does not constitute a loss of generality since it can always be achieved by the 
canonical construction given in \cite[Appendix]{lr}.
We assume $\Hom(\iota,\iota)=\7C\,\id_\iota$, i.e.\ the unit object $\iota$ is 
irreducible. 

In this work we will frequently deal with subcategories $\2S\subset\2C$. All such
subcategories will be assumed replete full. (A subcategory $\2S\subset\2C$ is
full iff $\Hom_\2S(\rho,\sigma)=\Hom_\2C(\rho,\sigma) \ \ \forall\rho,\sigma\in\2S$,
thus it is  determined by $\obj\,\2S$. A subcategory is replete iff 
$\rho\in\2S$ implies $\sigma\in\2S$ for all $\sigma\in\2C$ isomorphic to $\rho$.)
The replete full subcategories of $\2C$ form a lattice under inclusion, where
$\2S_1\subset\2S_2$ means $\obj\,\2S_1\subset\obj\,\2S_2$. 

Let now $\2S\subset\2C$ be a (replete full) symmetric subcategory closed under 
conjugates, direct sums and subobjects, the standard example being $\2D(\2C)$ by Prop.\ 
\ref{dc}. We do not assume $\2S\subset\2D(\2C)$ but we require that $\2S$ is \ul{even}
and refer to Subsect.\ \ref{super} for the supergroup case.
By the duality theorem of Doplicher and Roberts we have a unique compact group $G$ and an
invertible functor $F: \2S\rightarrow U(G)$. Here $U(G)$ is a category of finite 
dimensional continuous unitary 
representations of $G$, which is closed under subrepresentations and direct sums and 
which contains members of each isomorphism class of irreducible representations.
(Note that we did not specify the cardinalities of isomorphism classes in $U(G)$, since
they depend on the cardinalities in the given category $\2S$!)
The identity object of the category $U(G)$, viz.\ the space $\2H_0\cong\7C$ on which 
the trivial representation of $G$ `acts', contains a unit vector $\Omega$ such that the 
following identifications hold:
\be \Omega\boxtimes\psi=\psi\boxtimes\Omega=\psi\quad \forall \2H\in\obj\,U(G), \,
  \psi\in\2H.\ee

In order to avoid confusion with a later use of $\otimes$, the tensor product of objects 
in $F(\2S)=U(G)$, which are Hilbert spaces, will be denoted by $\boxtimes$ (as already
done above) and the product of objects $\rho,\sigma$ in $\2C$ by simple juxtaposition 
$\rho\sigma$. The composition and tensor product of arrows will be denoted by $\circ$ 
and $\times$, respectively, in both categories.
Let $\hat{G}$ be the set of all isomorphism classes of irreducible objects in $\2S$
or, equivalently by the duality theorem, of irreducible representations of $G$.
Let $\{\gamma_k, k\in\hat{G}\}$ be a section of objects in $\2S$ such that 
$\gamma_0=\iota$ and let $\2H_k=F(\gamma_k)$ be the images under the functor $F$. 
For every triple $k,l,m\in\hat{G}$ we choose an orthonormal basis 
\be \{V_{k,l}^{m,\alpha},\ \alpha=1,\ldots, N_{k,l}^m\} \ee
in the space $\Hom(\gamma_m,\gamma_k\gamma_l)$. (The latter space of arrows is in fact
a Hilbert space, but it should not be confused with the spaces $\2H_k, k\in\hat{G}$.)

The set $\hat{G}$ has an involution $k\mapsto\ol{k}$ which associates to every 
isomorphism class of representations of $G$ the conjugate class. By the isomorphism
between $\2S\cong U(G)$ this implies for our chosen section that $\gamma_{\ol{k}}$ is
conjugate to $\gamma_k$. Thus there are intertwiners 
$R_k\in\Hom(\iota,\gamma_{\ol{k}}\gamma_k),\ \ol{R}_k\in \Hom(\iota,\gamma_k\gamma_{\ol{k}})$
such that
\be \ol{R}_k^*\times\id_{\gamma_k}\mcirc\id_{\gamma_k}\times R_k=  \id_{\gamma_k}, 
  \quad\quad  R_k^*\times\id_{\gamma_{\ol{k}}}\mcirc\id_{\gamma_{\ol{k}}}\times
  \ol{R}_k=  \id_{\gamma_{\ol{k}}}. \label{conj-eq}\ee
Since this symmetric under $k\leftrightarrow\ol{k}, R\leftrightarrow\ol{R}$ one can 
choose $R_{\ol{k}}=\ol{R}_k, \ol{R}_{\ol{k}}=R_k$ for conjugate pairs of 
non-selfconjugate objects. For selfconjugate objects it is known that one can achieve
either $\ol{R}_k=R_k$ or $\ol{R}_k=-R_k$ depending on whether $\gamma_k$ is real or
pseudo-real. The above choices will be assumed in the sequel.

Now we define a new category $\2C\rtimes_0\2S$ in terms of the data 
$\2C, \2S, F, \hat{G}$.
\bdefin The category $\2C\rtimes_0\2S$ has the same objects as $\2C$ with the
same tensor product. The arrows in $\2C\rtimes_0\2S$ are defined by
\be \Hom_{\2C\rtimes_0\2S}(\rho,\sigma)=\bigoplus_{k\in\hat{G}} \
   \Hom_\2C(\gamma_k\rho,\sigma) \bigotimes \2H_k  \label{arr}\ee
with the obvious complex vector space structure. In order to economize on brackets we 
declare the precedence of products to be  $\otimes > \times > \circ > \bigotimes$, where
$\otimes, \bigotimes$ are different symbols for the tensor product in (\ref{arr}). 

Let $k,l\in\hat{G},\ T\in\Hom(\gamma_l\rho,\sigma), S\in\Hom(\gamma_k\sigma,\delta)$ and 
$\psi_k\in\2H_k,\ \psi_l\in\2H_l$. Then the composition of arrows in $\2C\rtimes_0\2S$ 
is defined by
\bea \lefteqn{\Hom_{\2C\rtimes_0\2S}(\rho,\delta) \ni S\otimes\psi_k\,\circ\,T\otimes\psi_l 
   = \bigoplus_{m\in\hat{G}} \sum_{\alpha=1}^{N_{k,l}^m} } \label{d-circ}\\    &&
   S\,\circ \,\id_{\gamma_k}\times T \,\circ\,
   V_{k,l}^{m,\alpha}\times\id_\rho \ \bigotimes \ F(V_{k,l}^{m,\alpha})^* 
   (\psi_k\boxtimes\psi_l) \nn\label{compos}\eea
and linear extension. Here $F$ is the embedding functor, thus $F(V_{k,l}^{m,\alpha})^*$
is a partial isometry from $\2H_k\boxtimes\2H_l$ onto $\2H_m$. 

Let $k,l\in\hat{G},\ S\in\Hom(\gamma_k\rho_1,\sigma_1), T\in\Hom(\gamma_l\rho_2,\sigma_2)$ and 
$\psi_k\in\2H_k,\ \psi_l\in\2H_l$. Then the tensor product of arrows in $\2C\rtimes_0\2S$
is defined by 
\bea \lefteqn{\Hom_{\2C\rtimes_0\2S}(\rho_1\rho_2,\sigma_1\sigma_2) \ni 
   S\otimes\psi_k\,\times\,T\otimes\psi_l = \bigoplus_{m\in\hat{G}} 
   \sum_{\alpha=1}^{N_{k,l}^m} } \label{d-times}\\   &&
   S \times T \,\circ\, \id_{\gamma_k} \times
   \ve(\gamma_l,\rho_1)\times\id_{\rho_2} \,\circ\, 
   V_{k,l}^{m,\alpha}\times\id_{\rho_1\rho_2} \ \bigotimes \ 
  F(V_{k,l}^{m,\alpha})^* (\psi_k\boxtimes\psi_l)  .\label{tens}\nn\eea
Finally, the $*$-operation of $\2C\rtimes_0\2S$ on the arrows
$S\otimes\psi_k\in\Hom_{\2C\rtimes_0\2S}(\rho,\sigma)$ with 
$S\in\Hom(\gamma_k\rho,\sigma), \psi\in\2H_k$  is defined by 
\be (S\otimes\psi_k)^* = R^*_k\times\id_\rho \,\circ\,
  \id_{\gamma_{\ol{k}}}\times S^* \ \bigotimes \ 
   \langle \psi_k\boxtimes\cdot\,, F(\ol{R}_k)\,\Omega \rangle. \label{star1}\ee
\label{maindef}\edefin
\rems 1. Tangle diagrams corresponding to the first tensor factor (which lives in the 
category $\2C$) in the definitions of $\circ, \times$ and $*$  are given in Figs.\ 
\ref{fig1} and \ref{fig5}. \\
2. A different choice for the orthonormal bases 
$\{V_{k,l}^{m,\alpha},\ \alpha=1,\ldots, N_{k,l}^m\}$ in $\Hom(\gamma_m,\gamma_k\gamma_l)$ 
does not affect the definition of $\circ, \times$, since the unitary matrices 
effecting the base change drop out. \\
3. The left tensor factor of (\ref{star1}) is in 
$\Hom(\gamma_{\ol{k}}\sigma,\rho)$, and $F(\ol{R}_{\gamma_k})\,\Omega$ 
is in $\2H_k\boxtimes\2H_{\ol{k}}$ such that contraction with $\psi_k$ yields a vector
in $\2H_{\ol{k}}$. Thus the entire expression is in $\Hom_{\2C\rtimes_0\2S}(\sigma,\rho)$ 
as it must be. \\
4. For every pair $\rho,\sigma$ there is an embedding of $\Hom_\2C(\rho,\sigma)$ into
$\Hom_{\2C\rtimes_0\2S}(\rho,\sigma)$ via $S\mapsto S\otimes\Omega$. Looking at the 
definitions of $\circ,\times$ in $\2C\rtimes_0\2S$ it is obvious that this gives rise to
a faithful functor from $\2C$ to $\2C\rtimes_0\2S$, thus $\2C$ can and will be 
considered as a subcategory of $\2C\rtimes_0\2S$.
Arrows in $\2C\rtimes_0\2S$ will be denoted $\tilde{S}, \tilde{T},\ldots$, but 
often we do not distinguish between $S\in\Hom_\2C(\rho,\sigma)$ and 
$S\otimes\Omega\in\Hom_{\2C\rtimes_0\2S}(\rho,\sigma)$. \\
5. By Frobenius reciprocity we have
$\dim\Hom(\gamma_k\rho,\sigma)=\dim\Hom(\gamma_k,\sigma\ol{\rho})<\infty$ and only finitely
many $k\in\hat{G}$ contribute, thus $\Hom_{\2C\rtimes_0\2S}(\rho,\sigma)$ is finite
dimensional. As a consequence of $\Hom(\gamma_k,\iota)=\{0\}$ for $k\ne e$ we obtain
\be \Hom_{\2C\rtimes_0\2S}(\iota,\iota)=\Hom(\iota,\iota)=\7C\,\id_\iota. \ee
6. A special case of (\ref{arr}) is 
\be \Hom_{\2C\rtimes_0\2S}(\iota,\gamma_k)=\Hom(\gamma_k,\gamma_k)\bigotimes\2H_k \ee
for $\gamma_k\in\2S$. Since the dimension $d_k\in\7N$ of $\gamma_k$ equals the dimension
of $\2H_k=F(\gamma_k)$, this implies $\gamma_k\cong d_k \iota$ in $\2C\rtimes_0\2S$. 
Thus $\gamma_k$ `disappears without a trace' in $\2C\rtimes_0\2S$ as far as the
irreducible objects are concerned. Furthermore, the spaces $\2H_k$ and 
$\Hom_{\2C\rtimes_0\2S}(\iota,\gamma_k)$
can be identified via $\psi\mapsto\id_{\gamma_k}\otimes\psi$. This allows us to consider
$\psi_k\in\2H_k$ also as a morphism in $\Hom_{\2C\rtimes_0\2S}(\iota,\gamma_k)$, which leads
to notational simplification. With $S\in\Hom(\gamma_k\rho,\sigma), \psi_k\in\2H_k$ it is an
easy consequence of (\ref{compos}) that 
$\Hom_{\2C\rtimes_0\2S}(\rho,\sigma)\ni S\otimes\psi_k=S\otimes\Omega\ \circ\ \id_{\gamma_k}\otimes\psi_k\,\times\,\id_\rho\otimes\Omega$. 
With the above identifications this
can also be written as $S\mcirc \psi_k\times\id_\rho$. In a sense, the new morphisms 
$\psi_k\in\Hom_{\2C\rtimes_0\2S}(\iota,\gamma_k)$ are the crucial point of Defin.\ 
\ref{maindef} and (\ref{arr}) simply reflects the fact that arrows can be composed.
It must of course still be proved that Defin.\ \ref{maindef} yields a BT$C^*$. \\
7. If $\2S\not\subset\2D$ then $\ve(\gamma,\rho)\circ\ve(\rho,\gamma)\ne\id_{\rho\gamma}$
for some $\gamma\in\2S, \rho\in\2C$. Thus there is another possible definition of 
$\times$ in $\2C\rtimes_0\2S$, replacing $\ve(\gamma_l,\rho_1)$ by
$\ve(\rho_1,\gamma_l)^{-1}$ in (\ref{tens}). For $\2S\subset\2D$ these definitions 
coincide. \\
8. Finally, we remark that there are similarities between our definition of 
$\2C\rtimes_0\2S$ and a construction \cite{r1} of a field algebra in algebraic quantum
field theory which preceded \cite{dr2} but where the main result of \cite{dr1} was 
assumed.

\begin{figure}
\[ \ba{ccc} \begin{tangle}
\step[1.5]\object{\delta} \\
\hh\step\hstep\id \\
\step\frabox{S} \\
\hh\step\id\step\id\step[.3]\mobj{\sigma} \\
\hh\step\id\step\frabox{T} \\
\hh\step[.2]\mobj{\gamma_k}\step[.8]\id\step[.3]\mobj{\gamma_l}\step[.7]\id\step\id \\
\step\Frabox{V_{k,l}^{m,\alpha}}\step\id \\
\hh\step[1.5]\id\step[1.5]\id \\
\step[1.5]\object{\gamma_m}\step[1.5]\object{\rho}
\end{tangle} 
& \quad\quad &
\begin{tangle}
\step[1.5]\object{\sigma_1}\step[3]\object{\sigma_2} \\
\hh\step[1.5]\id\step[3]\id \\
\step\frabox{S}\Step\frabox{T} \\
\step[.2]\mobj{\gamma_k}\step[.8]\id\hstep\mobj{\gamma_l}\hstep\xx\step\id \\
\step\Frabox{V_{k,l}^{m,\alpha}}\Step\id\step\id \\
\hh\step[1.5]\id\step[2.5]\id\step\id \\
\step[1.5]\object{\gamma_m}\step[2.5]\object{\rho_1}\step\object{\rho_2}
\end{tangle} \\
\hbox{Composition} && \hbox{Tensor Product} \ea \]
\caption{Composition and Tensor product of arrows in $\2C\rtimes_0\2S$}
\label{fig1}\end{figure}

\begin{figure}
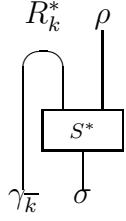

\[ \begin{tangle}
\hstep\object{R_k^*}\step[1.5]\object{\rho} \\
\hh\hcoev\step\id \\
\hh\id\step\id\step\id \\
\hh\id\step\frabox{S^*} \\
\hh\id\step[1.5]\id \\
\object{\gamma_{\ol{k}}}\step[1.5]\object{\sigma}
\end{tangle} \]
\caption{$*$-Operation on arrows}
\label{fig5}\end{figure}

\subsection{$C\rtimes_0\2S$ is a Tensor Category}
\blemma The operations $\circ, \times$ are bilinear and associative. \elemma
\prf Bilinearity is obvious. In order to prove associativity of $\circ$ consider
$S, T$ as in the definition (\ref{d-circ}) and $U\in\Hom(\gamma_n\eta,\rho)$. Then
\bea \lefteqn{(S\otimes\psi_k\mcirc T\otimes\psi_l)\mcirc U\otimes\psi_n =
  \bigoplus_{r\in\hat{G}} \sum_{m\in\hat{G}} \sum_{\alpha=1}^{N_{k,l}^m} 
   \sum_{\beta=1}^{N_{m,n}^r} } \nn\\ &&
   S\mcirc\id_{\gamma_k}\times T\mcirc\id_{\gamma_k\gamma_l}\times U \mcirc
   V_{k,l}^{m,\alpha}\times\id_{\gamma_n\eta}\mcirc V_{m,n}^{r,\beta}\times\id_{\eta} 
  \nn\\
 && \bigotimes\ (F(V_{m,n}^{r,\beta})^*\mcirc F(V_{k,l}^{m,\alpha})^*\times\id_{\2H_n})
  (\psi_k\boxtimes\psi_l\boxtimes\psi_n).  \eea
On the other hand
\bea \lefteqn{S\otimes\psi_k\mcirc (T\otimes\psi_l\mcirc U\otimes\psi_n) =
  \bigoplus_{r\in\hat{G}} \sum_{m\in\hat{G}} \sum_{\alpha=1}^{N_{l,n}^m} 
   \sum_{\beta=1}^{N_{k,m}^r} } \nn\\ &&
   S\mcirc\id_{\gamma_k}\times T\mcirc\id_{\gamma_k\gamma_l}\times U \mcirc
   \id_{\gamma_k}\times V_{l,n}^{m,\alpha}\times\id_\eta \mcirc V_{k,m}^{r,\beta}
  \times\id_\eta  \nn\\   && \bigotimes\ (F(V_{k,m}^{r,\beta})^* \mcirc \id_{\2H_k}\times
   F(V_{l,n}^{m,\alpha})^*) (\psi_k\boxtimes\psi_l\boxtimes\psi_n).  \eea
Since $F$ is a functor of $*$-categories we have
\bea F(V_{m,n}^{r,\beta})^*\mcirc F(V_{k,l}^{m,\alpha})^*\times\id_{\2H_n} &=& 
  F(V_{k,l}^{m,\alpha}\times\id_{\gamma_n}\mcirc V_{m,n}^{r,\beta})^*,\\
  F(V_{k,m}^{r,\beta})^* \mcirc \id_{\2H_k}\times F(V_{l,n}^{m,\alpha})^* &=&
  F(\id_{\gamma_k}\times V_{l,n}^{m,\alpha}\mcirc V_{k,m}^{r,\beta})^*. \eea
Since both
\be \{ V_{k,l}^{m,\alpha}\times\id_{\gamma_n}\mcirc V_{m,n}^{r,\beta},\ m\in\hat{G}, 
  \alpha=1,\ldots,N_{k,l}^m, \beta=1,\ldots, N_{m,n}^r \} \ee
and 
\be \{ \id_{\gamma_k}\times V_{l,n}^{m,\alpha} \mcirc V_{k,m}^{r,\beta},\ m\in\hat{G}, 
  \alpha=1,\ldots,N_{l,n}^m, \beta=1,\ldots,N_{k,m}^r \} \ee
are orthogonal bases of $\Hom(\gamma_r,\gamma_k\gamma_l\gamma_n)$, the two 
expressions coincide.
 
\begin{figure}
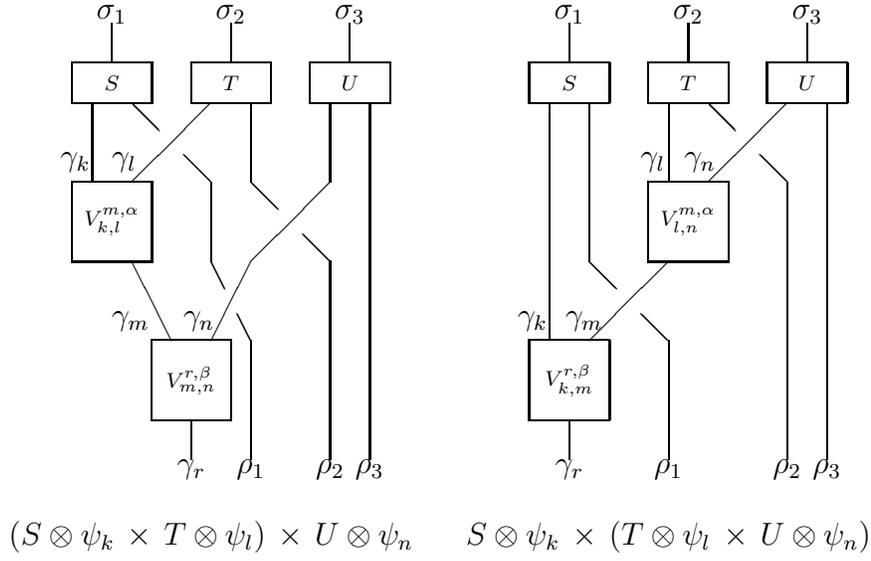

\[ \ba{ccc} \begin{tangle}
\step[1.5]\object{\sigma_1}\step[3]\object{\sigma_2}\step[3]\object{\sigma_3} \\
\hh\step[1.5]\id\step[3]\id\step[3]\id \\
\step\frabox{S}\Step\frabox{T}\Step\frabox{U} \\
\step[.2]\mobj{\gamma_k}\step[.8]\id\hstep\mobj{\gamma_l}\hstep\xx\step\id\Step\id\step\id \\
\step\Frabox{V_{k,l}^{m,\alpha}}\Step\id\step\xx\step\id \\
\step[1.5]\mobj{\gamma_m}\step[.5]\d\step[.3]\mobj{\gamma_n}\step[.7]\hxx\Step\id\step\id \\
\step[3]\Frabox{V_{m,n}^{r,\beta}}\step \id\Step\id\step\id \\
\hh\step[3.5]\id\step[1.5]\id\Step\id\step\id \\
\step[3.5]\object{\gamma_r}\step[1.5]\object{\rho_1}\Step\object{\rho_2}\step\object{\rho_3}
\end{tangle} &&
\begin{tangle} 
\step[1.5]\object{\sigma_1}\step[3]\object{\sigma_2}\step[3]\object{\sigma_3} \\
\hh\step[1.5]\id\step[3]\id\step[3]\id \\
\step\frabox{S}\Step\frabox{T}\Step\frabox{U} \\
\step\id\step\id\step[1.3]\mobj{\gamma_l}\step[.7]\id\step[.4]\mobj{\gamma_n}\step[.6]\xx\step\id \\
\step\id\step\id\Step\Frabox{V_{l,n}^{m,\alpha}}\Step\id\step\id \\
\step[.2]\mobj{\gamma_k}\step[.8]\id\step[.4]\mobj{\gamma_m}\step[.6]\xx\step[3]\id\step\id \\
\step\Frabox{V_{k,m}^{r,\beta}}\Step\id\step[3]\id\step\id \\
\hh\step[1.5]\id\step[2.5]\id\step[3]\id\step\id \\
\step[1.5]\object{\gamma_r}\step[2.5]\object{\rho_1}\step[3]\object{\rho_2}\step\object{\rho_3}
\end{tangle} \\ \\
(S\otimes\psi_k\,\times\,T\otimes\psi_l)\,\times\, U\otimes\psi_n &&
S\otimes\psi_k\,\times\,(T\otimes\psi_l\,\times\, U\otimes\psi_n) 
\ea \]
\caption{Associativity of $\times$}
\label{fig2}\end{figure}

The proof of associativity of $\times$ is similar. Let $S,T$ be as in (\ref{d-times}) 
and let $U\in\Hom(\gamma_n\rho_3,\sigma_3)$. Since writing down (and reading!) the formulae
would be rather tedious we express the parts of the summands which live in $\2C$
graphically, cf.\ Fig.\ \ref{fig2}. Thus
\bea \lefteqn{(S\otimes\psi_k\,\times\,T\otimes\psi_l)\,\times\, U\otimes\psi_n =
  \bigoplus_{r\in\hat{G}} \sum_{m\in\hat{G}} \sum_{\alpha=1}^{N_{k,l}^m} 
   \sum_{\beta=1}^{N_{m,n}^r}} \nn\\ 
   && (\mbox{Fig.\ \ref{fig2}, l.h.s.})\ \bigotimes \
  F(V_{k,l}^{m,\alpha}\times\id_{\gamma_n}\mcirc V_{m,n}^{r,\beta})^*
  (\psi_k\boxtimes\psi_l\boxtimes\psi_n).  \eea
On the other hand
\bea \lefteqn{S\otimes\psi_k\,\times\,(T\otimes\psi_l\,\times\, U\otimes\psi_n) =
  \bigoplus_{r\in\hat{G}} \sum_{m\in\hat{G}} \sum_{\alpha=1}^{N_{l,n}^m} 
 \sum_{\beta=1}^{N_{k,m}^r}} \nn\\ 
  && (\mbox{Fig.\ \ref{fig2}, r.h.s.})\ \bigotimes \
  F(\id_{\gamma_k}\times V_{l,n}^{m,\alpha}\mcirc V_{k,m}^{r,\beta})^* 
  (\psi_k\boxtimes\psi_l\boxtimes\psi_n).  \eea
By naturality the arrow $V_{l,n}^{m,\alpha}$ in the r.h.s.\ of Fig.\ \ref{fig2} 
can be pulled through the braiding, and the identity of the two expressions follows by 
the same argument as for $\circ$. \qed

\begin{figure}
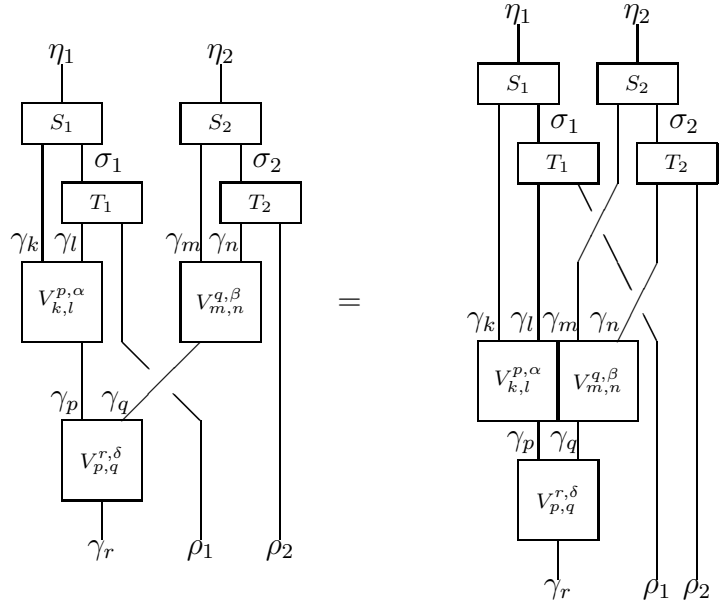

\[ \ba{ccc}  
\begin{tangle}
\step[1.5]\object{\eta_1}    \step[4]\object{\eta_2} \\
\hh\step[1.5]\id             \step[4]\id \\
\step\frabox{S_1}            \step[3]\frabox{S_2} \\
\hh\step\id\step\id\step[.3]\mobj{\sigma_1}\step[2.7] \id\step\id\step[.3]\mobj{\sigma_2} \\
\hh\step\id\step\frabox{T_1} \step[2]\id\step\frabox{T_2} \\
\hh\step[.2]\mobj{\gamma_k}\step[.8]\id\step[.3]\mobj{\gamma_l}\step[.7]\id\step\id
 \step\step[.1]\mobj{\gamma_m}\step[.9]\id\step[.2]\mobj{\gamma_n}\step[.8]\id\step\id \\
\step\Frabox{V_{k,l}^{p,\alpha}}\step\id \step[2]\Frabox{V_{m,n}^{q,\beta}}\step\id \\
\step[1.2]\mobj{\gamma_p}\step[.8]\id\step[.5]\mobj{\gamma_q}\step[.5]\xx                    \Step \id \\
\Step\Frabox{V_{p,q}^{r,\delta}}     \Step\id\Step\id \\
\hh\step[2.5]\id\step[2.5]\id\Step\id \\
\step[2.5]\object{\gamma_r}\step[2.5]\object{\rho_1}\Step\object{\rho_2} 
\end{tangle}
 & \quad = \quad &
\begin{tangle} 
\step[2.5]\object{\eta_1}\step[3]\object{\eta_2} \\
\hh\step[2.5]\id\step[3]\id \\
\hh\Step\frabox{S_1}\Step\frabox{S_2} \\
\hh\Step\id\step\id\step[.3]\mobj{\sigma_1}\step[1.7]\id\step\id\step[.3]\mobj{\sigma_2} \\
\hh\Step\id\step\frabox{T_1}\step\id\step\frabox{T_2} \\
\Step\id\step\id\step\hxx\step\id\step\id \\
\step\step[.2]\mobj{\gamma_k}\step[.8]\id\step[.3]\mobj{\gamma_l}\step[.7]\id\step[.1]\mobj{\gamma_m}\step[.9]\id\step[.3]\mobj{\gamma_n}\step[.7]\hxx\step\id \\
\Step\Frabox{V_{k,l}^{p,\alpha}}\step\Frabox{V_{m,n}^{q,\beta}}\step\id\step\id \\
\hh\Step\step[.2]\mobj{\gamma_p}\step[.8]\id\step[.3]\mobj{\gamma_q}\step[.7]\id\Step\id\step\id \\
\step[3]\Frabox{V_{p,q}^{r,\delta}}\Step\id\step\id \\
\hh\step[3.5]\id\step[2.5]\id\step\id \\
\step[3.5]\object{\gamma_r}\step[2.5]\object{\rho_1}\step\object{\rho_2}
\end{tangle}
\ea \]
\caption{$(S_1\otimes\psi_k\mcirc T_1\otimes\psi_l)\ \times\ (S_2\otimes\psi_m\mcirc T_2\otimes\psi_n)$}
\label{fig3}\end{figure}

\begin{figure}
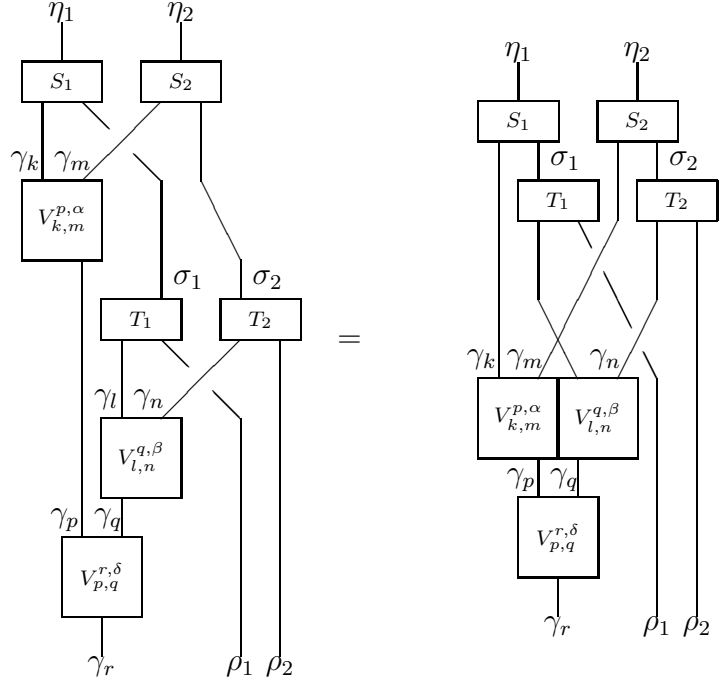

\[ \ba{ccc}  
\begin{tangle} 
\step[1.5]\object{\eta_1}\step[3]\object{\eta_2} \\
\hh\step[1.5]\id\step[3]\id \\
\step\frabox{S_1}\Step\frabox{S_2} \\
\step[.2]\mobj{\gamma_k}\step[.8]\id\step[.3]\mobj{\gamma_m}\step[.7]\xx\step\id \\
\step\Frabox{V_{k,m}^{p,\alpha}}\Step\id\step\d \\
\hh\Step\id\Step\id\step[.3]\mobj{\sigma_1}\step[1.7]\id\step[.3]\mobj{\sigma_2} \\
\hh\Step\id\step\frabox{T_1}\Step\frabox{T_2} \\
\Step\id\step[.3]\mobj{\gamma_l}\step[.7]\id\step[.3]\mobj{\gamma_n}\step[.7]\xx\step\id \\
\Step\id\step\Frabox{V_{l,n}^{q,\beta}}\Step\id\step\id \\
\hh\step\step[.2]\mobj{\gamma_p}\step[.8]\id\step[.3]\mobj{\gamma_q}\step[.7]\id\step[3]\id\step\id   \\
\Step\Frabox{V_{p,q}^{r,\delta}}\step[3]\id\step\id \\
\hh\step[2.5]\id\step[3.5]\id\step\id \\
\step[2.5]\object{\gamma_r}\step[3.5]\object{\rho_1}\step\object{\rho_2}
\end{tangle} 
 & \quad = \quad &
\begin{tangle} 
\step[2.5]\object{\eta_1}\step[3]\object{\eta_2} \\
\hh\step[2.5]\id\step[3]\id \\
\hh\Step\frabox{S_1}\Step\frabox{S_2} \\
\hh\Step\id\step\id\step[.3]\mobj{\sigma_1}\step[1.7]\id\step\id\step[.3]\mobj{\sigma_2} \\
\hh\Step\id\step\frabox{T_1}\step\id\step\frabox{T_2} \\
\Step\id\step\id\step\hxx\step\id\step\id \\
\step\step[.2]\mobj{\gamma_k}\step[.8]\id\step[.2]\mobj{\gamma_m}\step[.8]
\hbx(1,2){\put(0,2){\line(1,-2){1}}\put(0,0){\line(1,2){1}}}
\step[.3]\mobj{\gamma_n}\step[.7]\hxx\step\id \\
\Step\Frabox{V_{k,m}^{p,\alpha}}\step\Frabox{V_{l,n}^{q,\beta}}\step\id\step\id \\
\hh\Step\step[.2]\mobj{\gamma_p}\step[.8]\id\step[.3]\mobj{\gamma_q}\step[.7]\id\Step\id\step\id \\
\step[3]\Frabox{V_{p,q}^{r,\delta}}\Step\id\step\id \\
\hh\step[3.5]\id\step[2.5]\id\step\id \\
\step[3.5]\object{\gamma_r}\step[2.5]\object{\rho_1}\step\object{\rho_2}
\end{tangle}
\ea \]
\caption{$S_1\otimes\psi_k\,\times\,S_2\otimes\psi_m\ \mcirc \
   T_1\otimes\psi_l\,\times\,T_2\otimes\psi_n$}
\label{fig4}\end{figure}

\blemma The operations $\circ,\times$ satisfy the interchange law
\be (\tilde{S}_1\circ\tilde{T}_1)\times(\tilde{S}_2\circ\tilde{T}_2)=
  \tilde{S}_1\times\tilde{S}_2\,\circ\,\tilde{T}_1\times\tilde{T}_2, \label{intlaw}\ee
whenever the left hand side is defined. \elemma
\prf We compute
\bea \lefteqn{(S_1\otimes\psi_k\mcirc T_1\otimes\psi_l)\ \times\ (S_2\otimes\psi_m\mcirc
   T_2\otimes\psi_n) =  \bigoplus_{r\in\hat{G}} \sum_{p,q\in\hat{G}} 
  \sum_{\alpha=1}^{N_{k,l}^p} \sum_{\beta=1}^{N_{m,n}^q} 
  \sum_{\delta=1}^{N_{p,q}^r} } \nn\\ && (\mbox{Fig.\ \ref{fig3}, l.h.s.}) \
 \bigotimes\ F(V_{k,l}^{p,\alpha}\times V_{m,n}^{q,\beta} \mcirc V_{p,q}^{r,\delta})^* 
  (\psi_k\boxtimes\psi_l\boxtimes\psi_m\boxtimes\psi_n) \label{i1}\eea
and
\bea \lefteqn{S_1\otimes\psi_k\,\times\,S_2\otimes\psi_m\ \circ\ T_1\otimes\psi_l\,
  \times\,T_2\otimes\psi_n = \bigoplus_{r\in\hat{G}} \sum_{p,q\in\hat{G}} 
   \sum_{\alpha=1}^{N_{k,m}^p} \sum_{\beta=1}^{N_{l,n}^q}  
   \sum_{\delta=1}^{N_{p,q}^r} } \nn\\ && (\mbox{Fig.\ \ref{fig4}, l.h.s.}) \
 \bigotimes\ F(V_{k,m}^{p,\alpha}\times V_{l,n}^{q,\beta}\mcirc V_{p,q}^{r,\delta})^*
  (\psi_k\boxtimes\psi_m\boxtimes\psi_l\boxtimes\psi_n). \label{i2}\eea
(Since $\2S$ is a symmetric category we have used the symmetric braiding symbol for 
$\ve(\gamma_m,\gamma_l)$ in Fig.\ \ref{fig4}. We do not do this for braidings of 
$\gamma$'s with objects not in $\2S$ since we do not assume $\2S$ to be degenerate.)
By standard manipulations the left hand sides of Figs.\ \ref{fig3}, \ref{fig4} can be 
seen to equal the respective right hand sides. Next we transform (\ref{i2}) using
\bea \lefteqn{F(V_{k,m}^{p,\alpha}\times V_{l,n}^{q,\beta}\mcirc V_{p,q}^{r,\delta})^*
  (\psi_k\boxtimes\psi_m\boxtimes\psi_l\boxtimes\psi_n) = } \\ &&
   F(\id_{\gamma_k}\times\ve(\gamma_m,\gamma_l)\times\id_{\gamma_n}\mcirc 
   V_{k,m}^{p,\alpha}\times V_{l,n}^{q,\beta}\mcirc V_{p,q}^{r,\delta})^*
  (\psi_k\boxtimes\psi_l\boxtimes\psi_m\boxtimes\psi_n), \nn\eea
and observing that 
$\{ V_{k,l}^{p,\alpha}\times V_{m,n}^{q,\beta} \mcirc V_{p,q}^{r,\delta} \}$ and
$\{ \id_{\gamma_k}\times\ve(\gamma_m,\gamma_l)\times\id_{\gamma_n}\mcirc 
   V_{k,m}^{p,\alpha}\times V_{l,n}^{q,\beta}\mcirc V_{p,q}^{r,\delta} \}$
are orthonormal bases of $\Hom(\gamma_r,\gamma_k\gamma_l\gamma_m\gamma_n)$ 
(with $p,q\in\hat{G}$ and $\alpha, \beta,\delta$ in the obvious ranges) we are done. \qed

\blemma $\2C\rtimes_0\2S$ has conjugates and direct sums. \elemma
\prf Since the objects of $\2C\rtimes_0\2S$ are just those of $\2C$ the existence of
conjugates in $\2C\rtimes_0\2S$ follows from
\be R_k\in\Hom(\iota,\gamma_{\ol{k}}\gamma_k)\subset
   \Hom_{\2C\rtimes_0\2S}(\iota,\gamma_{\ol{k}}\gamma_k) \ee
and the fact that the conjugate equations clearly hold in $\2C\rtimes_0\2S$, too. In
the same way one shows that $\2C\rtimes_0\2S$ has direct sums. \qed

\subsection{The $*$-Operation}
\blemma The $*$-operation is antilinear and involutive. \elemma
\prf Antilinarity is obvious by definition. As to involutivity consider
$\tilde{S}=S\otimes\psi_k$ with $S\in\Hom(\gamma_k\rho,\sigma)$, $\psi\in\2H_k$.
Twofold application of the $*$-operation (\ref{star1}) yields
\be (S\otimes\psi_k)^{**} =
 R^*_{\ol{k}}\times\id_\sigma\mcirc\id_{\gamma_k\gamma_{\ol{k}}}\times S \mcirc 
  \id_{\gamma_k}\times R_k \times\id_\rho \ \bigotimes \ \langle \ol{\psi}_{\ol{k}}
  \boxtimes\cdot\,, F(\ol{R}_{\ol{k}})\,\Omega\rangle,\label{2star}\ee
where 
\be \ol{\psi}_{\ol{k}}=\langle \psi_k\boxtimes\cdot\,, F(\ol{R}_k)\,\Omega\rangle. \ee
The first tensor factor of (\ref{2star}) (which lives in $\2C$) can be transformed as 
follows:
\be \ba{ccccc} \begin{tangle}
\hstep\object{R^*_{\ol{k}}}\step[2.5]\object{\sigma} \\
\hh\hcoev\Step\id \\
\hh\id\step\id\step\frabox{S^{**}} \\
\hh\id\step\hev\step\id \\
\object{\gamma_k}\step[1.5]\object{R_k}\step[1.5]\object{\rho}
\end{tangle} 
& \quad=\quad &
\begin{tangle} 
\step[3]\object{\sigma} \\
\hh\step[3]\id \\
\Step\frabox{S} \\
\hh\hcoev\step\id\step\id \\
\hh\id\step\hev\step\id \\
\object{\gamma_k}\step[3]\object{\rho}
\end{tangle} 
& \quad=\ \ \pm \ \  &
\begin{tangle} 
\step\object{\sigma} \\
\hh\step\id \\
\frabox{S} \\
\hh\id\step\id \\
\object{\gamma_k}\step\object{\rho}
\end{tangle} 
\ea \ee
In the first step we have used interchange law and in the second step the first 
conjugate equation (\ref{conj-eq}). The possible appearance of the minus sign is due the
fact that $R^*_{\ol{k}}$ appears in (\ref{2star}) instead of $\ol{R}^*_k$.
In view of our choice of $R_{\ol{k}}=\pm\ol{R}_k$ the minus sign appears iff $k$ is 
selfconjugate and pseudoreal. Abbreviating the second factor in (\ref{2star}) (which
lives in $U(G)$) by $\ol{\ol{\psi}}_k$ we have
\be \langle a, \ol{\ol{\psi}}_k \rangle = \langle \ol{\psi}_{\ol{k}}\boxtimes a, 
   F(\ol{R}_{\ol{k}})\,\Omega \rangle \quad\forall a\in\2H_k. \ee
Inserting
\be \langle \ol{\psi}_{\ol{k}}, b\rangle =\langle F(\ol{R}_k)\,\Omega, \psi_k\boxtimes b
   \rangle \quad\forall b\in\2H_{\ol{k}}\ee
we have
\bea \langle a, \ol{\ol{\psi}}_k \rangle &=& \langle F(\ol{R}_k)\,\Omega\,\boxtimes\,a, 
   \psi_k\,\boxtimes\,F(\ol{R}_{\ol{k}})\,\Omega\rangle \nn\\
  &=& \langle \Omega\boxtimes a, F(\ol{R}^*_k\times\id_{\gamma_k}\mcirc
   \id_{\gamma_k}\times\ol{R}_{\ol{k}}) \,\psi_k\boxtimes\Omega \rangle. \eea
Now 
$\ol{R}^*_k\times\id_{\gamma_k}\mcirc\id_{\gamma_k}\times\ol{R}_{\ol{k}}=\pm \id_{\gamma_k}$, 
thus $F(\ldots)=\pm\id_{\2H_k}$. With $\Omega\boxtimes a=a$
and $\psi_k\boxtimes\Omega=\psi_k$ we have 
$\langle a, \ol{\ol{\psi}}_k \rangle = \pm\langle a,\psi_k\rangle$ and therefore 
$\ol{\ol{\psi}}_k= \pm \psi_k$. Also here the minus sign appears iff
$k$ is selfconjugate and pseudoreal. In any case the two minus signs cancel and 
we obtain $(S\otimes\psi_k)^{**}=(S\otimes\psi_k)$. \qed

\begin{figure}
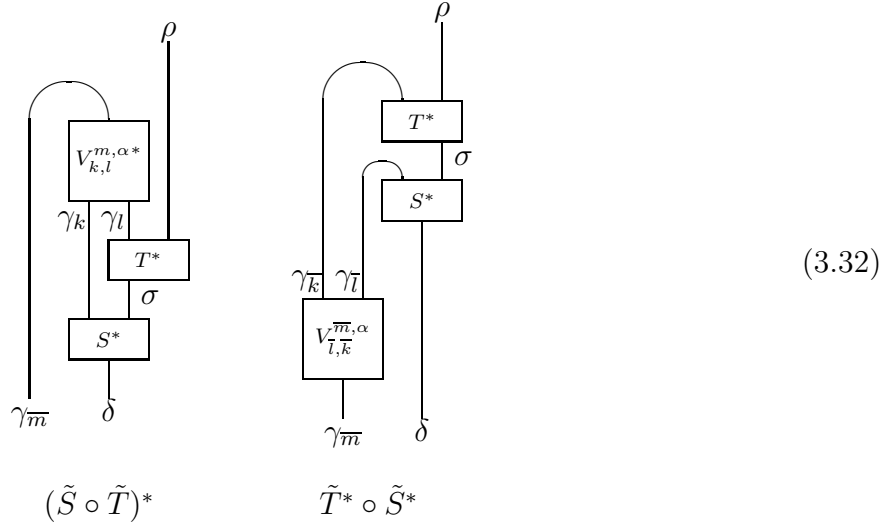

\be \ba{ccc} \begin{tangle}
\step[3.5]\object{\rho} \\
\coev\step[1.5]\id \\
\id\step[1.5]\Frabox{{V_{k,l}^{m,\alpha}}^*}\step\id \\
\hh\id\hstep\step[.2]\mobj{\gamma_k}\step[.8]\id\step[.3]\mobj{\gamma_l}\step[.7]\id\step\id \\
\hh\id\step[1.5]\id\step\frabox{T^*} \\
\hh\id\step[1.5]\id\step\id\step[.3]\mobj{\sigma} \\
\hh\id\step[1.5]\frabox{S^*} \\
\hh\id\Step\id \\
\object{\gamma_{\ol{m}}}\Step\object{\delta}
\end{tangle} 
 & \enspace\enspace\enspace\enspace &
\begin{tangle} 
\step[4]\object{\rho} \\
\step\coev\step\id \\
\hh\step\id\Step\frabox{T^*} \\
\hh\step\id\step\hcoev\step\id\step[.3]\mobj{\sigma} \\
\hh\step\id\step\id\step\frabox{S^*} \\
\step[.2]\mobj{\gamma_{\ol{k}}}\step[.8]\id\step[.3]\mobj{\gamma_{\ol{l}}}\step[.7]\id\step[1.5]\id \\
\step\Frabox{V_{\ol{l},\ol{k}}^{\ol{m},\alpha}}\step[1.5]\id \\
\hh\step[1.5]\id\Step\id \\
\step[1.5]\object{\gamma_{\ol{m}}}\Step\object{\delta}
\end{tangle}  \\ \\
(\tilde{S}\circ\tilde{T})^* && \tilde{T}^*\circ\tilde{S}^*
\ea \ee
\caption{Compatibility of $*$ and $\circ$}
\label{fig7}\end{figure}

\blemma The $*$-operation is contravariant, i.e.\ 
$(\tilde{S}\circ\tilde{T})^*=\tilde{T}^*\circ\tilde{S}^*$ whenever the left hand side 
is defined.  \elemma
\prf Let $S\in\Hom(\gamma_k\sigma,\delta), T\in\Hom(\gamma_l\rho,\sigma)$ and 
$\psi_k\in\2H_k, \psi_l\in\2H_l$ and apply the $*$-operation (\ref{star1}) to 
$\tilde{S}\circ\tilde{T}=S\otimes\psi_k\mcirc T\otimes\psi_l$
as defined by (\ref{compos}). We obtain
\bea (\tilde{S}\circ\tilde{T})^* &=& \bigoplus_{m\in\hat{G}} \sum_{\alpha=1}^{N_{k,l}^m}\
   R^*_m\times\id_\rho\mcirc\id_{\gamma_{\ol{m}}}\times{V_{k,l}^{m,\alpha}}^*\times
  \id_\rho \mcirc \id_{\gamma_{\ol{m}}\gamma_k}\times T^* \mcirc \id_{\gamma_{\ol{m}}}
 \times S^* \nn\\ && \bigotimes \ \langle F(V_{k,l}^{m,\alpha})^* (\psi_k\boxtimes\psi_l)
  \,\boxtimes\, \cdot\,, F(\ol{R}_m)\,\Omega \rangle. \label{s1}\eea
On the other hand,
\bea \tilde{T}^*\circ\tilde{S}^* &=& \bigoplus_{m\in\hat{G}} 
  \sum_{\alpha=1}^{N_{\ol{k},\ol{l}}^{\ol{m}}} \
   R^*_l\times\id_{\rho} \mcirc\id_{\gamma_{\ol{l}}}\times T^* \mcirc
  \id_{\gamma_{\ol{l}}}\times R^*_k\times\id_\sigma \mcirc 
  V_{\ol{l},\ol{k}}^{\ol{m},\alpha} \times S^*   \nn\\ && \bigotimes \
  F(V_{\ol{l},\ol{k}}^{\ol{m},\alpha})^* (\ol{\psi}_{\ol{l}}\boxtimes
   \ol{\psi}_{\ol{k}}), \label{s2}\eea
where 
\be \ol{\psi}_{\ol{l}}= \langle \psi_l\boxtimes\cdot\,, F(\ol{R}_l)\,\Omega\rangle \ee
and similarly for $\ol{\psi}_{\ol{k}}$. The left tensor factors of (\ref{s1}) and 
(\ref{s2}) are represented in Fig.\ \ref{fig7}.

As to the right hand factors of (\ref{s1}) and (\ref{s2}) which live in $\2H_{\ol{m}}$
and which we abbreviate $\psi_1, \psi_2$, respectively, we have for all 
$a\in\2H_{\ol{m}}$:
\bea \langle a, \psi_1\rangle &=& \langle F(V_{k,l}^{m,\alpha})^* (\psi_k\boxtimes\psi_l)
  \,\boxtimes\, a, F(\ol{R}_m)\,\Omega \rangle \nn\\
  &=& \langle \psi_k\boxtimes\psi_l\boxtimes a, F(V_{k,l}^{m,\alpha}\times
  \id_{\gamma_{\ol{m}}}\mcirc \ol{R}_m)\,\Omega \rangle \label{psi1}\eea
and
\bea \langle a, \psi_2\rangle &=& \langle a, F(V_{\ol{l},\ol{k}}^{\ol{m},\alpha})^* 
  (\ol{\psi}_{\ol{l}}\boxtimes \ol{\psi}_{\ol{k}}) \rangle 
  = \langle F(V_{\ol{l},\ol{k}}^{\ol{m},\alpha}) a, 
  \ol{\psi}_{\ol{l}}\boxtimes \ol{\psi}_{\ol{k}} \rangle \nn\\
  &=& \langle\psi_k\boxtimes \psi_l \boxtimes\,F(V_{\ol{l},\ol{k}}^{\ol{m},\alpha}) a,
  [F(\ol{R}_l)\Omega]_{23} \, [F(\ol{R}_k)\Omega]_{14} \rangle \label{psi2}\\
  &=& \langle\psi_k\boxtimes \psi_l \boxtimes\,F(V_{\ol{l},\ol{k}}^{\ol{m},\alpha}) a,
  F(\id_{\gamma_k}\times\ol{R}_l\times\id_{\gamma_{\ol{k}}}\mcirc \ol{R}_k)\, \Omega
   \rangle \nn\\
  &=& \langle\psi_k\boxtimes \psi_l \boxtimes a, 
   F(\id_{\gamma_k\gamma_l}\times V_{\ol{l},\ol{k}}^{\ol{m},\alpha *} \mcirc
  \id_{\gamma_k}\times\ol{R}_l\times\id_{\gamma_{\ol{k}}}\mcirc \ol{R}_k)\,\Omega\rangle.
\nn\eea
The fourth equality in (\ref{psi2}) follows from the following computation in
$\2H_k\boxtimes\2H_l\boxtimes\2H_{\ol{l}}\boxtimes\2H_{\ol{k}}$:
\bea \lefteqn{[F(\ol{R}_l)\Omega]_{23} \, [F(\ol{R}_k)\Omega]_{14} =
   \sigma_{12}\circ\sigma_{23}(F(\ol{R}_l)\Omega\,\boxtimes\, F(\ol{R}_k)\Omega)=} \\
 && F(\ve(\gamma_l\gamma_{\ol{l}},\gamma_k)\times\id_{\gamma_{\ol{k}}}\mcirc
   \ol{R}_l \times \ol{R}_k)\,\Omega =
  F(\id_{\gamma_k}\times\ol{R}_l\times\id_{\gamma_{\ol{k}}}\mcirc \ol{R}_k)\,\Omega,
\nn\eea
where in the last step we have used the interchange law.

Now we observe that 
$\{ W_{k,l}^{m,\beta}, \beta=1,\ldots,N_{k,l}^m \}$
with
\be W_{k,l}^{m,\beta}= \ol{R}_m^* \mcirc
   \id_{\gamma_m}\times V_{\ol{l},\ol{k}}^{\ol{m},\beta *} \mcirc
   \id_{\gamma_m\gamma_{\ol{l}}}\times R_k\times\id_{\gamma_l} \mcirc 
   \id_{\gamma_m}\times R_l \ee
is an orthonormal basis in $\Hom(\gamma_m,\gamma_k\gamma_l)$. Since the choice of such a
basis is irrelevant we can replace $V_{k,l}^{m,\alpha}$ in (\ref{s1}) by 
$W_{k,l}^{m,\alpha}$. Using the conjugate equations (\ref{conj-eq}) one then 
easily verifies that (\ref{s1}) and (\ref{s2}) coincide. \qed

\begin{figure}
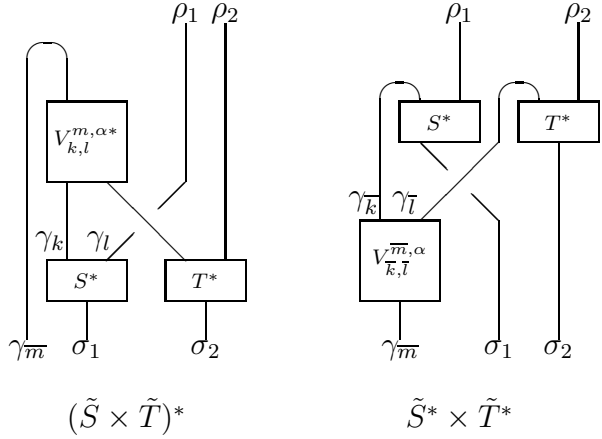

\[ \ba{ccc}  
\begin{tangle} 
\step[4]\object{\rho_1}\step\object{\rho_2} \\
\hh\hcoev\step[3]\id\step\id \\
\hh\id\step\id\step[3]\id\step\id \\
\id\step\Frabox{V_{k,l}^{m,\alpha *}}\Step\id\step\id \\
\id\step[.2]\mobj{\gamma_k}\step[.8]\id\step[.5]\mobj{\gamma_l}\step[.5]\x\step\id \\
\hh\id\step\frabox{S^*}\Step\frabox{T^*} \\
\hh\id\step[1.5]\id\step[3]\id \\
\object{\gamma_{\ol{m}}}\step[1.5]\object{\sigma_1}\step[3]\object{\sigma_2} 
\end{tangle} 
 & \enspace\enspace\enspace\enspace &
\begin{tangle} 
\step[3]\object{\rho_1}\step[3]\object{\rho_2} \\
\step\hcoev\step\id\step\hcoev\step\id \\
\hh\step\id\step\frabox{S^*}\step\id\step\frabox{T^*} \\
\step[.2]\mobj{\gamma_{\ol{k}}}\step[.8]\id\step[.3]\mobj{\gamma_{\ol{l}}}\step[.7]\xx\step[1.5]\id \\
\step\Frabox{V_{\ol{k},\ol{l}}^{\ol{m},\alpha}}\Step\id\step[1.5]\id \\
\hh\step[1.5]\id\step[2.5]\id\step[1.5]\id \\
\step[1.5]\object{\gamma_{\ol{m}}}\step[2.5]\object{\sigma_1}\step[1.5]\object{\sigma_2} 
\end{tangle} \\ \\
  (\tilde{S}\times\tilde{T})^* && \tilde{S}^*\times\tilde{T}^*
\ea \]
\caption{Compatibility of $*$ and $\times$}
\label{fig6}\end{figure}

\blemma The $*$-operation is monoidal, i.e.\ $(S\times T)^*=S^*\times T^*$. \elemma
\prf Let $S\in\Hom(\gamma_k\rho_1,\sigma_1), T\in\Hom(\gamma_l\rho_2,\sigma_2), \psi_k\in\2H_k, \psi_l\in\2H_l$. 
Then
\be (\tilde{S}\times\tilde{T})^*=\bigoplus_{m\in\hat{G}} \sum_{\alpha=1}^{N_{k,l}^m} \
  (\mbox{Fig.\ \ref{fig6}, l.h.s.}) \ \bigotimes \ 
  \langle F(V_{k,l}^{m,\alpha})^* (\psi_k\boxtimes\psi_l)
  \,\boxtimes\, \cdot\,, F(\ol{R}_m)\,\Omega \rangle. \label{s3}\ee
On the other hand,
\be \tilde{S}^*\times\tilde{T}^*=\bigoplus_{m\in\hat{G}} 
  \sum_{\alpha=1}^{N_{\ol{k},\ol{l}}^{\ol{m}}} \
  (\mbox{Fig.\ \ref{fig6}, r.h.s.}) \ \bigotimes \
  F(V_{\ol{k},\ol{l}}^{\ol{m},\alpha})^* (\ol{\psi}_{\ol{k}}\boxtimes
   \ol{\psi}_{\ol{l}}). \label{s4}\ee

Using the interchange law several times, the right hand side of Fig.\ \ref{fig6} is 
shown to equal
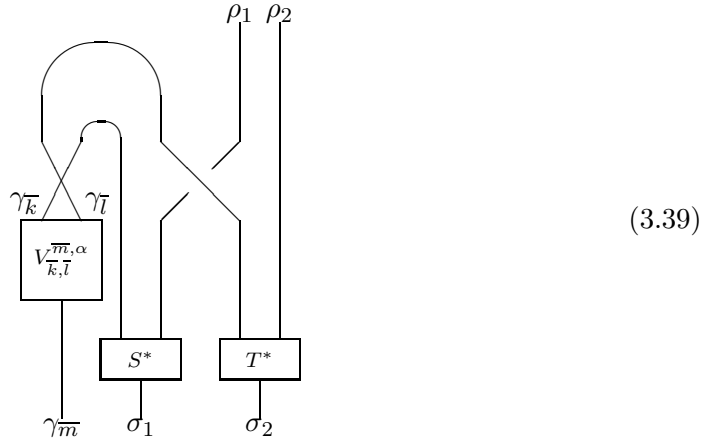
\be \begin{tangle}
\step[6]\object{\rho_1}\step\object{\rho_2} \\
\step\mcoev\Step\id\step\id \\
\hh\step\id\step\hcoev\step\id\Step\id\step\id \\
\step[.2]\mobj{\gamma_{\ol{k}}}\step[.8]\hbx(1,2){\put(0,2){\line(1,-2){1}}\put(0,0){\line(1,2){1}}} \step[.1]\mobj{\gamma_{\ol{l}}}\step[.9]\id\step\x\step\id \\
\step\Frabox{V_{\ol{k},\ol{l}}^{\ol{m},\alpha}}\step\id\step\id\Step\id\step\id \\
\hh\step[1.5]\id\step[1.5]\id\step\id\Step\id\step\id \\
\hh\step[1.5]\id\step[1.5]\frabox{S^*}\Step\frabox{T^*} \\
\hh\step[1.5]\id\Step\id\step[3]\id \\
\step[1.5]\object{\gamma_{\ol{m}}}\Step\object{\sigma_1}\step[3]\object{\sigma_2} 
\end{tangle}
\ee
which differs from the left hand side of Fig.\ \ref{fig6} only by a replacement of the
basis 
\be \{ R_m^* \mcirc \id_{\gamma_{\ol{m}}}\times V_{k,l}^{m,\alpha *}, \ \alpha=1,\ldots,
   N_{k,l}^m \} \ee
of $\Hom(\gamma_{\ol{m}}\gamma_k\gamma_l,\iota)$ by
\be \{ R_l^*\mcirc\id_{\gamma_{\ol{l}}}\times R_k^*\times\id_{\gamma_l} \mcirc
   (\ve(\gamma_{\ol{k}},\gamma_{\ol{l}})\circ V_{\ol{k},\ol{l}}^{\ol{m},\alpha})\times
   \id_{\gamma_k\gamma_l},\  \alpha=1,\ldots, N_{k,l}^m \}. \label{b2}\ee
Concerning the right hand sides the calculation proceeds as in the preceding lemma. 
The only difference is that in (\ref{s4}) 
$F(V_{\ol{k},\ol{l}}^{\ol{m},\alpha})^*(\ol{\psi}_{\ol{k}}\boxtimes\ol{\psi}_{\ol{l}})$
appears in contrast to 
$F(V_{\ol{l},\ol{k}}^{\ol{m},\alpha})^*(\ol{\psi}_{\ol{l}}\boxtimes\ol{\psi}_{\ol{k}})$
in (\ref{s2}). But this is compensated for by the $\ve(\gamma_{\ol{k}},\gamma_{\ol{l}})$
in (\ref{b2}). \qed

\blemma The $*$-operation of $\2C\rtimes_0\2S$ is positive. 
Thus $\2C\rtimes_0\2S$ is a $C^*$-tensor category. \elemma
\prf Let $\tilde{S}\in\Hom_{\2C\rtimes_0\2S}(\rho,\sigma)$. It is sufficient to prove that 
the vanishing of $(\tilde{S}^*\tilde{S})_e$, i.e.\ the component  in (\ref{arr}) with 
$k=e$ (the $G$-invariant part, see below), implies $\tilde{S}=0$. Let thus 
\be S=\bigoplus_{k\in\hat{G}} \sum_i \ S_k^i\bigotimes \psi_k^i,\quad
   S_k^i\in\Hom(\gamma_k\rho,\sigma),\ \psi_k\in\2H_k. \ee
(We must sum over an index $i$ in order to allow for elements of 
$\Hom_{\2C\rtimes_0\2S}(\rho,\sigma)$ which are not simple tensors.) Then 
\[  (\tilde{S}^*\tilde{S})_e=\sum_{k,l\in\hat{G}} \sum_{i,j} R_k^*\times
  \id_\rho \mcirc \id_{\gamma_{\ol{k}}}\times {S_k^i}^* \mcirc \id_{\gamma_{\ol{k}}}
  \times S_l^j \mcirc V_{\ol{k},l}^e \times\id_\rho \]
\be \bigotimes \ F(V_{\ol{k},l}^e)^* 
   (\langle\psi_k^i\boxtimes\cdot, F(\ol{R}_k)\,\Omega\rangle  \boxtimes \psi_l^j). \ee
Now, the space $\Hom(\gamma_e,\gamma_{\ol{k}}\gamma_l)=\Hom(\iota,\gamma_{\ol{k}}\gamma_l)$ is
one dimensional for $l=k$ and trivial otherwise. Since the choice of an orthonormal
basis in this space does not matter we can choose 
$V^e_{\ol{k},k}=d(k)^{-1/2} \,R_k$. Here the numerical factor involving the dimension
$d(k)=d(\gamma_k)>0$ \cite{lr} is necessary in order for $V$ to be isometric. Then
\[  (\tilde{S}^*\tilde{S})_e=\sum_{k\in\hat{G}} \frac{1}{d(k)} \sum_{i,j} 
  R_k^*\times \id_\rho \mcirc \id_{\gamma_{\ol{k}}}\times ({S_k^i}^* \circ S_k^j) \mcirc 
  R_k \times\id_\rho \]
\be \bigotimes \ F(R_k)^* 
   (\langle\psi_k^i\boxtimes\cdot, F(\ol{R}_k)\,\Omega\rangle \boxtimes \psi_k^j). \ee
Considering the $\Hom(\rho,\rho)$-valued bilinear form on $\Hom(\gamma_k\rho,\rho)$ 
\be (S, T)\mapsto\langle S,T \rangle_k\ = R_k^*\times\id_\rho \mcirc 
  \id_{\gamma_{\ol{k}}}\times (S^* \circ T) \mcirc R_k \times\id_\rho, \ee
positivity of the $*$-operation of $\2C$ implies that $\langle S,S \rangle_k=0$ iff
$\id_{\gamma_{\ol{k}}}\times S \mcirc R_k \times\id_\rho=0$. By Frobenius reciprocity
this is the case iff $S=0$, thus $\langle\cdot,\cdot\rangle_k$ is positive definite. 
Furthermore,
\bea F(R_k)^*(\langle\psi_k^i\boxtimes\cdot, F(\ol{R}_k)\,\Omega\rangle \boxtimes\psi_k^j)
 &=& \langle \Omega, F(R_k)^*(\langle\psi_k^i\boxtimes\cdot, F(\ol{R}_k)\,
    \Omega\rangle \boxtimes\psi_k^j) \rangle \ \Omega\nn\\
  &=& \langle F(R_k)\,\Omega, \langle\psi_k^i\boxtimes\cdot, F(\ol{R}_k)\,\Omega
    \rangle \boxtimes\psi_k^j \rangle \ \Omega\nn\\
  &=& \langle \psi_k^i\boxtimes F(R_k)\,\Omega, F(\ol{R}_k)\,\Omega\boxtimes\psi_k^j
    \rangle \ \Omega \nn\\
  &=& \langle \psi_k^i\boxtimes\Omega, F(\id_{\gamma_k}\times R_k^*\mcirc \ol{R}_k\times
   \id_{\gamma_k})\,\Omega\boxtimes\psi_k^j \rangle \ \Omega\nn\\
  &=& \langle \psi_k^i,\psi_k^j \rangle_{\2H_k} \ \Omega, \eea
where we have used the conjugate equations. Thus also
\be (\tilde{S}^*\tilde{S})_e  =  \sum_{k\in\hat{G}} \frac{1}{d(k)} \sum_{i,j} 
  \langle \psi_k^i, \psi_k^j \rangle_{\2H_k}\ \langle S_k^i, S_k^j\rangle_k\,\bigotimes\,
  \Omega \ee
is positive definite since it is the sum of the tensor product of such maps, 
and $\tilde{S}^*\tilde{S}$ vanishes iff $\tilde{S}=0$. The second claim follows by 
Prop.\ \ref{cstar}. \qed

Summing up we have proved
\bprop $\2C\rtimes_0\2S$ is a $C^*$-tensor category with conjugates and direct
sums. \eprop
\rem If $\2S\subset\2D$ we can consider the crossed product $\2D\rtimes_0\2S$, 
which is a full subcategory of $\2C\rtimes_0\2D$. It is interesting to note that 
$\2D\rtimes_0\2S$ can be defined also if $\2S\not\subset\2D$, namely as the full 
subcategory of $\2C\rtimes_0\2S$ whose objects are those in $\2D$. It is obvious
that for $\2S\subset\2D$ this notation is consistent with the crossed product in the
sense of Defin.\ \ref{maindef}. Thus also for $\2S\not\subset\2D$ we obtain a 
$C^*$-tensor category $\2D\rtimes_0\2S$ with conjugates and direct sums. It turns out,
however, that we obtain nothing new in this way. For, by Frobenius reciprocity in
$C^*$-tensor categories \cite{lr} we have 
$\dim\Hom_\2C(\gamma_k\rho,\sigma)=\dim\Hom_\2C(\gamma_k,\sigma\ol{\rho})$
In view of $\sigma\ol{\rho}\in\2S$ we have $\Hom_\2C(\gamma_k\rho,\sigma)=\{0\}$
whenever $\gamma_k\not\in\2S$. Thus the direct sum in (\ref{arr}) effectively runs only
over the $k$ such that $\gamma_k\in\2D$, which implies 
$\2D\rtimes_0\2S=\2D\rtimes_0(\2D\cap\2S)$. 
Therefore we are left with the crossed product of a symmetric tensor category by a full
subcategory.

\subsection{Braidings, Subobjects and Uniqueness}
\blemma The braiding $\ve$ of $\2C$ lifts to a braiding for $\2C\rtimes_0\2S$ iff 
$\2S\subset\2D$. \label{braid}\elemma
\prf Define 
$\tilde{\ve}(\rho,\sigma)=\ve(\rho,\sigma)\otimes\Omega\in\Hom_{\2C\rtimes_0\2S}(\rho\sigma,\sigma\rho)$.
That $\tilde{\ve}$ satisfies the relations 
\bea \tilde{\ve}(\rho,\sigma_1\sigma_2) &=& \id_{\sigma_1}\times\tilde{\ve}
   (\rho,\sigma_2)\mcirc \tilde{\ve}(\rho,\sigma_1)\times\id_{\sigma_2}, \\
   \tilde{\ve}(\rho_1\rho_2,\sigma)&=&\tilde{\ve}(\rho_1,\sigma)\times\id_{\rho_2}\mcirc
   \id_{\rho_1}\times\tilde{\ve}(\rho_2,\sigma) \eea
is obvious since these relations hold in $\2C$.
It remains to show that $\tilde{\ve}$ is natural w.r.t.\ both variables 
also in the extended category. Assuming $\2S\subset\2D$ we will prove 
\be \tilde{S}\times\id_\rho \mcirc \tilde{\ve}(\rho,\sigma) =
   \tilde{\ve}(\rho,\eta) \mcirc \id_\rho\times\tilde{S} \label{nat1}\ee 
in $\2C\rtimes_0\2S$ with $\tilde{S}\in\Hom_{\2C\rtimes_0\2S}(\sigma,\eta)$. The proof of
naturality w.r.t.\ the other variable is similar, and the general result follows by
the interchange law (\ref{intlaw}). Now, in more explicit terms the left hand side of 
(\ref{nat1}) amounts to (with $S\in\Hom(\gamma_k\sigma,\eta)$)
\bea S\otimes\psi_k \,\times\,\id_\rho\otimes\Omega \ \circ\ \ve(\rho,\sigma)\otimes
  \Omega &=& (S\times\id_\rho)\,\otimes\,\psi_k \ \circ\ \ve(\rho,\sigma)\otimes \Omega
 \nn\\ 
  &=& S\times\id_\rho \mcirc \id_{\gamma_k}\times\ve(\rho,\sigma)\ \bigotimes\ \psi_k
\eea
and the right hand side to
\bea \ve(\rho,\eta)\otimes\Omega \ \circ \ \id_\rho\otimes\Omega \times S\otimes\psi_k
  &=& \ve(\rho,\eta)\otimes\Omega \ \circ \ [\id_\rho\times S \mcirc \ve(\gamma_k,\rho)
   \times\id_\sigma ] \otimes\psi_k \nn\\
  &=& \ve(\rho,\eta)\mcirc \id_\rho\times S\mcirc \ve(\gamma_k,\rho)\times\id_\sigma 
   \ \bigotimes \ \psi_k. \eea
That these expressions coincide is seen by the following calculation for the $\2C$-parts.

\be \ba{ccccc} \begin{tangle}
\hstep\object{\eta}\step[1.5]\object{\rho} \\
\hh\hstep\id\step[1.5]\id \\
\hh\frabox{S}\step\id \\
\hh\id\step\id\step\id \\
\id\step\hxx \\
\object{\gamma_k}\step\object{\rho}\step\object{\sigma}
\end{tangle} 
& \quad=\quad &
\begin{tangle} 
\hstep\object{\eta}\step[1.5]\object{\rho} \\
\hh\hstep\id\step[1.5]\id \\
\hh\frabox{S}\step\id \\
\id\step\hxx \\
\hxx\step\id\\
\hxx\step\id\\
\object{\gamma_k}\step\object{\rho}\step\object{\sigma}
\end{tangle} 
& \quad=\quad &
\begin{tangle} 
\object{\eta}\step\object{\rho} \\
\hxx \\
\hh\id\step\id \\
\hh\id\step\frabox{S} \\
\hh\id\step\id\step\id \\
\hxx\step\id \\
\object{\gamma_k}\step\object{\rho}\step\object{\sigma}
\end{tangle} 
\ea \ee
In the second step of this computation we have used the naturality of the braiding
in $\2C$, and the first step is legitimate if $\ve_M(\gamma_k,\rho)=\id_{\gamma_k\rho}$.
This holds for all $\rho\in\2C$ if $\2S\subset\2D$ since then all $\gamma_k$ are 
degenerate.
Now assume $\2S\not\subset\2D$, i.e.\ there is a $\gamma_k\in\2S$ which has non-trivial
monodromy with some $\rho\in\2C$. Let now $\eta\prec\gamma_k\sigma$ and 
$S\in\Hom(\gamma_k\sigma,\eta)$. Reversing the above argument we see that naturality of
the braiding $\tilde{\ve}(\rho,\sigma)$ in $\2C\rtimes_0\2S$ fails for 
$\tilde{S}=S\otimes\psi_k\in\Hom_{\2C\rtimes_0\2S}(\sigma,\eta)$. \qed\\
\rem It is instructive to relate this result to what happens in the quantum field 
framework
\cite{khr2,mue4}. There the observables $\2A$ are extended by fields implementing the 
sectors in a symmetric semigroup $\Delta$ of DHR endomorphisms and the localized sectors
of $\2A$ are extended to the fields $\2F$. If $\Delta$ contains non-degenerate sectors 
then the extension $\tilde{\rho}$ of at least one sector $\rho$ is solitonic, i.e.\ 
localized only in a half-space. But it is well known that for solitons there is no 
braiding.

As observed in remark 6 after Defin.\ \ref{maindef}, the objects $\gamma_k\in\2S$
decompose into multiples of $\iota$ in $\2C\rtimes_0\2S$. But in $\2C\rtimes_0\2S$ also 
other irreducible objects $\rho\in\2C$ may become reducible in the sense that 
$\Hom_{\2C\rtimes_0\2S}(\rho,\rho)\supsetneq\7C\,\id_\rho$. In this
case the subobjects are not already present in $\2C$. Thus $\2C\rtimes_0\2S$ 
will in general not be closed under subobjects. 
There is a canonical procedure \cite[Appendix]{lr}, yielding for every 2-category $\2C$ 
a 2-category $\ol{\2C}$ which is closed under subobjects and contains $\2C$ as a 
full subcategory. Since we are concerned only with the special (and
more familiar) case of tensor categories, we give a fairly explicit description below.

\bdefin The closure $\ol{\2C}$ of a tensor category $\2C$ w.r.t.\ subobjects has as
objects pairs $(\rho,E)$ where $\rho\in\obj\,\2C$ and $E=E^2=E^*\in\Hom_\2C(\rho,\rho)$. 
The morphisms in $\ol{\2C}$ are given by 
\be \Hom_{\ol{\2C}}((\rho,E),(\sigma,F)) = \{T\in\Hom_\2C(\rho,\sigma)\ |\ 
   T=T\circ E=F\circ T \} = F\mcirc\Hom_\2C(\rho,\sigma)\mcirc E, \ee
and the composition of morphisms, where defined, is the one of $\2C$. The identity 
morphisms are given by $\id_{(\rho,E)}=E$. The tensor product is
$(\rho,E)\,(\sigma,F)=(\rho\sigma, E\times F)$ for the objects and the one
of $\2C$ for the morphisms. The embedding of $\2C$ in $\ol{\2C}$ is given by
$\rho\mapsto (\rho,\id_\rho)$ and the identity map on the arrows. \label{clossub}\edefin
\rem With this definition $(\rho,E)$ is a subobject of $\rho=(\rho,\id_\rho)$ in view of
$E\in\Hom((\rho,E),(\rho,\id_\rho))$ and $E\circ E^*=E, E^*\circ E=E=\id_{(\rho,E)}$. 
Assume a subobject $\rho_1\prec\rho$ exists in $\2C$ with $V\in\Hom_\2C(\rho_1,\rho)$ 
isometric. Then $\rho_1$ is isomorphic in $\ol{\2C}$ to $(\rho,E)$, where $E=V\circ V^*$.
Indeed, on one hand $V\in\Hom((\rho_1,\id_{\rho_1}),(\rho,E))$ since 
$V=V\circ\id_{\rho_1}=\id_{(\rho,E)}\circ V=E\circ V=V\circ V^*\circ V=V$. 
On the other hand, $V$ is unitary (in $\ol{\2C}$ !) since $V^*\circ V=\id_{\rho_1}$ 
and $V\circ V^*=E=\id_{(\rho,E)}$. 
If $\2C$ has conjugates then also $\ol{\2C}$ has conjugates. For, if 
$\rho, \ol{\rho}, R, \ol{R}$ satisfy the conjugate equations, then $(\ol{\rho},\ol{E})$
is a conjugate for $(\rho,E)$. Here 
\be \ol{E}=R^*\times\id_{\ol{\rho}}\mcirc\id_{\ol{\rho}}\times E\times\id_{\ol{\rho}}
   \mcirc\id_{\ol{\rho}}\times \ol{R} \label{olE}\ee
is easily verified to be an orthogonal projection in $(\ol{\rho},\ol{\rho})$, and
$R_{(\rho,E)}=\ol{E}\times E\mcirc R, \ol{R}_{(\rho,E)}=E\times\ol{E}\mcirc\ol{R}$
satisfy the conjugate equations. If $\2C$ is obtained from a subcategory $\2C_0$ by
adding morphisms, $\rho$ is irreducible in $\2C_0$ and $\ol{\rho},R, \ol{R}$ is a
solution of the conjugate equations in $\2C_0$ then with the above it is easy to see
that $\rho,\ol{\rho},R,\ol{R}$ is a standard solution in $\ol{\2C}$. Finally, given 
$V\in\Hom(\rho,\tau), W\in\Hom(\sigma,\tau)$ with $V\circ V^*+W\circ W^*=\id_\tau$ (thus 
$\tau\cong\rho\oplus\sigma$) and given projections $E\in\Hom(\rho,\rho), F\in\Hom(\sigma,\sigma)$
it is easy to verify that $(\tau, V\circ E\circ V^*+W\circ F\circ W^*)$ is a direct sum of 
$(\rho,E)$ and 
$(\sigma,F)$. Thus, if $\2C$ is closed w.r.t.\ subobjects then the obvious embedding
functor $\2C\rightarrow\ol{\2C}$ is essentially surjective. Since it is also full and
faithful, $\2C$ and $\ol{\2C}$ are equivalent as categories, cf.\ 
\cite[Sect.\ IV.4]{cwm}. That this is in fact an equivalence of tensor categories 
requires an additional argument for which we refer, e.g., to \cite{y}.

\bdefin $\2C\rtimes\2S=\ol{\2C\rtimes_0\2S}$. $\2C$ is identified with a
subcategory of $\2C\rtimes\2S$ via the embedding $\rho\mapsto(\rho,\id_\rho)$,
$\Hom(\rho,\sigma)\ni S\mapsto S\otimes\Omega\in\Hom((\rho,\id_\rho),(\sigma,\id_\sigma))$.
\edefin

\btheor $\2C\rtimes\2S$ is a $C^*$-tensor category with conjugates,
direct sums and subobjects. If $\2S\subset\2D$ then $\2C\rtimes\2S$ is braided. 
If $\2C$ is rational then so is $\2C\rtimes\2S$. \etheor
\prf As shown above, closing under subobjects does not affect the property of being 
closed under direct sums. Since an object $\rho$ has the same finite dimension in 
$\2C\rtimes\2S$ as in $\2C$, it decomposes into finitely many subobjects in 
$\2C\rtimes\2S$. Thus $\2C\rtimes\2S$ is rational if $\2C$ is.
It only remains to prove that the braiding of $\2C\rtimes_0\2S$ given by Lemma 
\ref{braid} if $\2S\subset\2D$ extends uniquely to the closure under subobjects. 
This was shown for symmetric tensor
categories in \cite{dr1} and works also in the braided case. We sketch the argument.
Consider $\rho,\sigma\in\obj\,\2C=\obj\,\2C\rtimes_0\2S$ and 
$E\in\Hom(\rho,\rho), F\in\Hom(\sigma,\sigma)$. Defining
\be \ve((\rho,E),(\sigma,F))=F\times E\mcirc\ve(\rho,\sigma)\mcirc E\times F,
\label{braid2}\ee
it is easily verified that we obtain a braiding for $\2C\rtimes\2S$ which satisfies
naturality w.r.t.\ both variables. \qed

\bprop Up to isomorphism of tensor categories, the category $\2C\rtimes\2S$ does not 
depend on the choice of the section $\{\gamma_l, l\in\hat{G}\}$ and of the functor $F$.
If $\2S\subset\2D$ then this isomorphism respects the braiding. \label{uniq}\eprop
\prf Let $\{\gamma_k, k\in\hat{G}\}, \{\gamma'_k, k\in\hat{G}\}$ be two sections of 
$\hat{G}$ in $\2S$ and let $F, F'$ be functors embedding $\2S$ into the category of 
Hilbert spaces. Denote the corresponding categories by 
$\2C\rtimes_0^{(\gamma,F)}\2S,\2C\rtimes_0^{(\gamma',F')}\2S$. We know that there are 
unitaries $W_k\in\Hom(\gamma_k,\gamma'_k)$ as well as a natural transformation 
$\{U_\rho:F(\rho)\rightarrow F'(\rho), \rho\in\hat{G}\}$ from $F$ to $F'$ 
with the $U_\rho$s being unitaries. Then the linear maps 
$\Hom_{\2C\rtimes_0^{(\gamma,F)}\2S}(\rho,\sigma)\rightarrow\Hom_{\2C\rtimes_0^{(\gamma',F')}\2S}(\rho,\sigma)$ defined by
\be S\otimes \psi_k \ \mapsto \ S\circ W_k\times\id_\rho \otimes U_k \psi_k,\quad
   S\in\Hom(\gamma_k\rho,\sigma), \psi_k\in\2H_k=F(\gamma_k) \ee
are isomorphisms. The easy proof that these maps define a (braided) tensor $*$-functor 
(obviously invertible)
from $\2C\rtimes_0^{(\gamma,F)}\2S$ to $\2C\rtimes_0^{(\gamma',F')}\2S$ which is the 
identity on the objects is left to the reader. Finally, isomorphic categories have 
isomorphic closures under subobjects. 
{}\vspace{1cm} \qed\\
\rem The functor $F$ is unique up to a natural transformation, the latter being in
one-to-one correspondence to the elements of $G$. The role of the compact group $G$ for 
the category $\2C\rtimes\2S$ we will thoroughly clarified in the next section.

\subsection{$G$-Symmetry}
By the DR duality theorem (or the Tannaka-Krein duality,
taking the existence of a representation functor $F$ for granted) the Hilbert spaces 
$\2H_k, k\in\hat{G}$ carry unitary representations $\pi_k(\cdot)$ of $G$. We define
an action of $G=\gal(\2S)$ on the morphisms of $\2C\rtimes_0\2S$ and thus of 
$\2C\rtimes\2S$ by
\be \alpha_g(S\otimes\psi_k)=S\otimes \pi_k(g)\psi_k, \quad S\in\Hom(\gamma_k\rho,\sigma). 
\label{alpha1}\ee
For the objects $(\rho,E)$ of $\2C\rtimes\2S=\ol{\2C\rtimes_0\2S}$ we define
\be \alpha_g((\rho,E))=(\rho,\alpha_g(E)), \label{alpha2}\ee
where $\alpha_g(E)$ is defined in (\ref{alpha1}). 

\bdefin Let $\2T\subset\2S$ be [B/S]T$C^*$s. Then $\aut_\2T(\2S)$ is the group of 
automorphisms (invertible [braided/symmetric] tensor $*$-endofunctors) of $\2S$
which leave $\2T$ stable. \edefin

\blemma The map $g\mapsto\alpha_g$ is a homomorphism of $G$ into 
$\aut_\2C(\2C\rtimes S)$. \elemma
\prf Using the definitions (\ref{compos}), (\ref{tens}) and the functoriality of $F$
one easily verifies that
\be \alpha_g(\tilde{S}\bullet\tilde{T})=\alpha_g(\tilde{S})\bullet\alpha_g(\tilde{T})
 \quad\mbox{where}\quad \bullet\in\{\circ,\times\}. \label{alpha0}\ee
In order to show that $\alpha_g$ is a functor it remains to show that 
$T\in\Hom_{C\rtimes\2S}((\rho,E),(\sigma,F))$ implies
\be \alpha_g(T)\in\Hom_{C\rtimes\2S}(\alpha_g((\rho,E)),\alpha_g((\sigma,F))). \ee
This is true due to $\alpha_g(T)\in\Hom_{C\rtimes_0\2S}(\rho,\sigma)$ and 
\be \alpha_g(T)=\alpha_g(T)\circ\alpha_g(E)=\alpha_g(F)\circ\alpha_g(T), \ee
where we have used (\ref{alpha0}). 
$\alpha_g$ is a tensor functor since (\ref{alpha0}) for $\bullet=\times$ implies
\be \alpha_g((\rho,E)\,(\sigma,F))=(\rho\sigma,\alpha_g(E\times F))=
  \alpha_g((\rho,E))\alpha_g((\sigma,F)). \ee
Finally, saying that $\alpha_g$ is a braided tensor functor is equivalent to the 
equation
\be \alpha_g(\ve((\rho,E),(\sigma,F)))=\ve(\alpha_g((\rho,E)),\alpha_g((\sigma,F))), 
\label{btf}\ee
which follows immediately from (\ref{braid2}) and the $G$-invariance of 
$\ve(\rho,\sigma)$. The homomorphism property of $g\mapsto\alpha_g$ is obvious and thus 
also the invertibility of $\alpha_g$. Clearly, $\alpha_g$ acts trivially on $\2C$. \qed

\bprop For every $\alpha\in\aut_\2C(\2C\rtimes S)$ there is $g\in G=\gal(\2S)$
such that $\alpha=\alpha_g$. Thus $\aut_\2C(\2C\rtimes S)\cong\gal(\2S)$.
\label{relgal}\eprop
\prf Let $\alpha\in\aut_\2C(\2C\rtimes S)$. Then $\alpha(\rho)=\rho$ for 
$\rho\in\2C$ implies $\alpha(T)\in\Hom_{\2C\rtimes S}(\rho,\sigma)$ if 
$T\in\Hom_{\2C\rtimes S}(\rho,\sigma)$. As before, we write $\tilde{T}=T\otimes\psi_k$
with $T\in\Hom_\2C(\gamma_k\rho,\sigma),\ \psi_k\in\2H_k$ also as 
$T\circ\psi_k\times\id_\rho$, where $\psi_k$ is interpreted as an element of
$\Hom_{\2C\rtimes S}(\iota,\gamma_k)$. Then 
$\alpha(\tilde{T})=T\circ\alpha(\psi_k)\times\id_\rho$ since $\alpha$ acts trivially on
the morphisms in $\2C$. Thus $\alpha$ is determined by the actions on the Hilbert spaces
$\2H_k$, which are clearly linear. Due to 
$\alpha(\psi^*\psi')=\psi^*\psi'\propto\id_\iota\in\2C$ for 
$\psi,\psi'\in\2H_k$ these actions are unitary, which then is true
for all spaces $\Hom_{\2C\rtimes S}(\iota,\gamma)$. If 
$\gamma,\gamma'\in\2S, V\in\Hom_\2C(\gamma,\gamma')$ and $\psi\in F(\gamma)$ then 
$\psi'=F(V)\psi\in F(\gamma')$ and $\alpha(\psi')=F(V)\alpha(\psi)$.
Thus $\alpha$ acts on the spaces $\Hom(\iota,\gamma), \gamma\in\2S$ like a natural 
transformation of the functor $F: \2S\rightarrow\2H$ to itself. 
But the latter are in one-to-one
correspondence to the elements of $G=\gal(\2S)$ \cite{dr1}. \qed

From now on we identify $G=\aut_\2C(\2C\rtimes S)$.
\bdefin Let $H\subset G$ be a subgroup. Then $(\2C\rtimes\2S)^H$ is the sub-tensor
category of $\2C\rtimes\2S$ consisting of $H$-invariant objects and morphisms. (That 
this really is a tensor category follows from the functoriality of $\alpha_g$.) \edefin

\blemma $(\2C\rtimes\2S)^G$ is equivalent to $\2C$. \label{CSG}\elemma
\prf A morphism 
$T\in\Hom_{\2C\rtimes\2S}((\rho,E),(\sigma,F))\subset \Hom_{\2C\rtimes_0\2S}(\rho,\sigma)$
is $G$-invariant iff it is in $\Hom_\2C(\rho,\sigma)$.
An object $(\rho,E)$ of $\2C\rtimes\2S$ is in $(\2C\rtimes\2S)^G$ iff $E$ is 
$G$-invariant iff $E\in\Hom_\2C(\rho,\rho)$. Thus $(\2C\rtimes\2S)^G$ is isomorphic to the
closure $\ol{\2C}$ of $\2C$ under subobjects. The latter is equivalent to $\2C$ since
$\2C$ is by assumption closed w.r.t.\ subobjects. (Recall the remark following Defin.\ 
\ref{clossub}.) \qed\\
\rems 1. The fact that $(\2C\rtimes\2S)^{\aut_\2C(\2C\rtimes S)}\simeq\2C$ justifies
calling $\2C\subset\2C\rtimes\2S$ a Galois extension of BT$C^*$s. This line of thought
will be continued in Subsect.\ \ref{galcorr}.\\
2. If $\rho\in\2C$ is irreducible then 
$\Hom_{\2C\rtimes\2S}(\rho,\rho)^G=\Hom_\2C(\rho,\rho)=\7C\,\id_\rho$, thus $G$ acts ergodically
on $\Hom_{\2C\rtimes\2S}(\rho,\rho)$. Now for irreducible $\rho\in\2C$ the obvious dimension
consideration 
\be \dim\Hom_\2C(\gamma_k\rho,\rho) \le d_k \quad\forall k\in\hat{G} \ee
together (\ref{arr}) implies that the irreducible representation $\pi_k$ of $G=\gal(\2S)$
occurs in $\Hom_{\2C\rtimes\2S}(\rho,\rho)$ with multiplicity at most $d_k$ (equivalently,
the corresponding spectral subspace has dimension at most $d_k^2$.) This is 
an instance of a well-known general result in the theory of ergodic compact group
actions on von Neumann algebras, cf.\ \cite[Prop.\ 2.1]{hls} or \cite[I]{wa}.

\sectreset{Galois Correspondence and the Modular Closure}\label{sect4}
Throughout the section $\2C$ is a BT$C^*$, $\2S\subset\2C$ is a ST$C^*$ and
$G=\aut_\2C(\2C\rtimes S)\cong\gal(\2S)$.
Having defined the semidirect product $\2C\rtimes\2S$ and established its uniqueness,
we will now prove some non-trivial properties. We continue to assume $\2S$ to be even
and will make clear which results require $\2S\subset\2D$.

\subsection{The Modular Closure $\2C\rtimes\2D$}
The following technical lemma can be distilled from \cite[I, Sect. 11]{wa}, but we
give the easy direct proof.
\blemma Let $N$ be a finite dimensional semisimple $\7C$-algebra and let 
$g\mapsto\alpha_g\in\aut\,N$ be an ergodic action of a group $G$. Then $N$ is 
isomorphic to the tensor product of its center $Z(N)$ with a full matrix algebra:
\be N \cong M_n \otimes Z(N) \cong \underbrace{M_n \oplus M_n \oplus \ldots 
   \oplus M_n}_{\mbox{$d$ terms}}, \ee
where $d=\dim\,Z(N)$ and $M_n$ denotes the simple algebra of complex $n\times n$ 
matrices. Let $E, F$ be minimal (i.e.\ one dimensional) projections in $N$. Then there 
is $g\in G$ such that $\alpha_g(E)\cong_N F$, i.e.\ there is $V\in N$ such that
$VV^*=F, V^*V=\alpha_g(E)$. \label{ergod}\elemma
\prf Let $E$ be a projection in $N$. Since $N$ is a von Neumann algebra it 
contains the projection $\ol{E}=\bigvee_{g\in G} \alpha_g(E)$, which clearly is
non-trivial and $G$-invariant. Therefore $\ol{E}\in N^G=\7C\11$ and thus 
$\ol{E}=\11$. Applying this to the (finite) set of minimal projections in $Z(N)$ we see
that $G$ acts transitively on the set of minimal central projections of $N$. Since the
dimension of such a projection is invariant under an automorphism of $N$, all simple
blocks of $N$ have the same rank.

Let $E, F$ be minimal projections in $N$ and let $\tilde{E},\tilde{F}$ be the (unique)
minimal projections in $Z(N)$ such that $E\le\tilde{E}, F\le\tilde{F}$. Then there is
$g\in G$ such that $\alpha_g(\tilde{E})=\tilde{F}$, thus $\alpha_g(E)\le\tilde{F}$.
This implies $\alpha_g(E)\cong_N F$ since all one dimensional projections in the factor
$\tilde{F}N$ are equivalent. \qed

\bprop Let $\rho\in\2C$ be irreducible. Then all irreducible subobjects $\rho_i$ of 
$\rho$ in $\2C\rtimes\2S$ occur with the same multiplicity and have the same dimension.
If $\2S\subset\2D$, thus $\2C\rtimes\2S$ is braided, then all $\rho_i$ have the same
twist as $\rho$, and they are either all degenerate or all non-degenerate according 
to whether $\rho$ is degenerate or non-degenerate. \label{crucial}\eprop
\prf $\Hom_{\2C\rtimes\2S}(\rho,\rho)$ is a finite dimensional von Neumann algebra and
$\Hom_{\2C\rtimes\2S}(\rho,\rho)^G=\Hom_\2C(\rho,\rho)=\7C\,\id_\rho$. Thus the lemma applies and
the first claim of the proposition follows from the result that all simple blocks of 
$\Hom_{\2C\rtimes\2S}(\rho,\rho)$ have the same rank. Let 
$E, F\in\Hom_{\2C\rtimes\2S}(\rho,\rho)$ be minimal projections corresponding to the 
irreducible subobjects $(\rho, E), (\rho, F)$ of $\rho$ and let $g, V$ be as in the 
lemma. Then $(\rho,\alpha_g(E))$ is equivalent to $(\rho,F)$ since $V$ is a unitary in 
$\Hom_{\2C\rtimes\2S}((\rho,\alpha_g(E)),(\rho,F))$. The dimension of $\rho$ being 
defined \cite{lr} via $d_\rho\,\id_\iota=R_\rho^*\circ R_\rho$ and $R_{(\rho,E)}$ being 
given as in the remark after Defin.\ \ref{clossub}, the independence of $d_{(\rho, E)}$ 
on $E$ follows from the transitivity of the $G$-action on the set of minimal central
projections. 

Assuming now $\2S\subset\2D$ it follows similarly that the twist is the 
same for all subobjects. If $\{ V_i\in\Hom(\rho_i,\rho) \}$ is a family of isometries such 
that $V_i^*\circ V_j=\delta_{i,j}\id_{\rho_i}$ and $\sum_i V_i\circ V_i^*=\id_\rho$ 
where the $\rho_i$ are irreducible in $\2C\rtimes\2S$, then 
$\kappa(\rho)=\sum_i V_i\circ\kappa(\rho_i)\circ V_i^*$. Since 
$\kappa(\rho_i)=\omega\,\id_{\rho_i}$ for some $\omega\in\7C$, this implies 
$\kappa(\rho)=\omega\id_\rho$ and thus $\omega(\rho)=\omega(\rho_i)\ \forall i$. 
Since $\alpha_g$ is a braided tensor functor (\ref{btf}), $(\rho, E)$ is degenerate iff
$(\rho,\alpha_g(E))\cong(\rho, F)$ is degenerate. Thus the subobjects $\rho_i$ are either
all degenerate or all non-degenerate. Since an object is degenerate iff
all subobjects are degenerate, cf.\ Prop.\ \ref{dc}, we conclude that the subobjects
are degenerate iff $\rho$ is degenerate. \qed\\
\rem That the decomposition of a degenerate object yields only degenerate objects was 
known before, cf.\ Prop.\ \ref{dc}, and for degenerate $\rho$ the result on the 
multiplicities and dimensions of the irreducible subobjects reduces to a well known 
result on group representations, as will be shown in the next subsection.
But for the non-degenerate objects, which have no group theoretic interpretation, the
above result in new and crucial for the rest of the paper. A detailed analysis of how
an irreducible non-degenerate object of $\2C$ decomposes in $\2C\rtimes\2S$ will be 
given in Subsect.\ \ref{abel} for the case where $\gal(\2S)$ is an abelian group.

\bcoro If $\2S\subset\2D$ then $\2D(\2C\rtimes\2S)\cong\2D\rtimes\2S$. 
\label{gcoro}\ecoro
\prf We have to show that starting from $\2C$ the operations of taking the crossed 
product with $\2S$ and of picking the full subcategory of degenerate objects commute.
Now we observe
\be \2D(\2C\rtimes\2S)=\2D(\ol{\2C\rtimes_0\2S})\cong
   \ol{\2D(\2C\rtimes_0\2S)}\cong\ol{\2D\rtimes_0\2S}=\2D\rtimes\2S, \ee
where the equalities hold by definition. The first isomorphism follows
since an irreducible subobject $(\rho,E)$ of $\rho$ is degenerate iff $\rho$ is 
degenerate, and the second isomorphism is true since 
$\2D(\2C\rtimes_0\2S)=\2D\rtimes_0\2S$. \qed

Even though further machinery will be developed below, we are already in a position to 
state one of our main results, which in fact provided the motivation for the entire
paper.

\btheor $\2C\rtimes\2D$ is non-degenerate. Thus every irreducible degenerate object 
is equivalent to $\iota$. If $\2C\rtimes\2D$ is rational (which follows if $\2C$ is
rational) then $\2C\rtimes\2D$ is modular. \label{main}\etheor
\prf By the proposition we have $\2D(\2C\rtimes\2D)\cong\2D\rtimes\2D$. Now,
all objects of $\2D\rtimes\2D$ are multiples of the identity, cf.\ 
remark 6 after Defin.\ \ref{maindef} Thus there are no irreducible degenerate objects 
in $\2C\rtimes\2D$ which are inequivalent to $\iota$. The rest follows from the 
discussion in Sect.\ 2. \qed

This result motivates the following
\bdefin The modular closure of a braided tensor $*$-category with conjugates, direct sums
and subobjects is $\ol{\ol{\2C}}=\2C\rtimes\2D$. \edefin
The terminology {\it closure} is justified by the fact that $\2D(\ol{\ol{\2C}})$ is 
trivial, which implies that the modular closure $\ol{\ol{\2C}}$ does not admit further 
crossed products (with braiding).

\subsection{Galois Correspondence}\label{galcorr}
Turning now to the study of categories $\2E$ sitting between $\2C$ and $\2C\rtimes\2S$ 
we begin with those of the form $(\2C\rtimes\2S)^H$ where $H\subset G$.
\blemma Let $H\subset G$ be a subgroup and let $\ol{H}$ be its closure in $G$. Then
$(\2C\rtimes\2S)^H=(\2C\rtimes\2S)^{\ol{H}}$ is a [B]T$C^*$. \label{int1}\elemma
\prf That the fixpoint categories under $H$ and $\ol{H}$ are the same follows from
continuity of $\pi_k$ in (\ref{alpha1}). If $E\in\Hom_{\2C\rtimes\2S}(\rho,\rho)$ 
is $H$-invariant then also $\ol{E}$ defined in (\ref{olE}) is $H$-invariant, thus 
$(\2C\rtimes\2S)^H$ has conjugates. That $(\2C\rtimes\2S)^H$ has direct sums and 
subobjects is seen similarly. In order to prove closedness of $(\2C\rtimes\2S)^H$
under the $*$-operation we have to show that $T^*$ is $H$-invariant if $T$ is.
In view of (\ref{star1}) we have
\be \alpha_g((S\otimes\psi_k)^*) =
   R^*_k\times\id_\rho \,\circ\, \id_{\gamma_{\ol{k}}}\times S^* \ \bigotimes \ 
  \pi_{\ol{k}}(g) \langle \psi_k\boxtimes\cdot\,, F(\ol{R}_k)\,\Omega \rangle. \ee
That $(S\otimes\psi_k)^*$ is $H$-invariant follows from the following calculation with 
$g\in H$ and $\psi_k\in\2H_k^H$:
\bea \lefteqn{\pi_{\ol{k}}(g) \langle \psi_k\boxtimes\cdot\,, F(\ol{R}_k)\,\Omega \rangle
   = \langle \psi_k\boxtimes\cdot\,, \pi_k(g)\times\pi_{\ol{k}}(g)\,F(\ol{R}_k)\,\Omega 
   \rangle} \nn\\
   && = \langle \psi_k\boxtimes\cdot\,, F(\ol{R}_k)\pi_0(g)\,\Omega \rangle 
  = \langle \psi_k\boxtimes\cdot\,, F(\ol{R}_k)\,\Omega \rangle. \eea
We have used that $\{\pi_k(g), k\in\hat{G} \}$ is a natural transformation of $F$
and that $\pi_0$ is the trivial representation. The restriction of the braiding of
$\2C\rtimes\2S$ to $(\2C\rtimes\2S)^H$ is, of course, a braiding. \qed

In order to prove that all T$C^*$s between $\2C$ and $\2C\rtimes\2S$ are of the form
$(\2C\rtimes\2S)^H$ we need the following
\blemma Let $\2E$ be a T$C^*$ such that $\2C\subset\2E\subset\2C\rtimes\2S$. 
With the identification of the spaces Hilbert $\2H_k=F(\gamma_k)$ and
$\Hom_{\2C\rtimes\2S}(\iota,\gamma_k)=\Hom(\gamma_k,\gamma_k)\otimes\2H_k$ via
$\psi_k\mapsto\id_{\gamma_k}\otimes\psi_k$ we have
\be \Hom_\2E(\rho,\sigma) = \bigoplus_{k\in\hat{G}} \ \Hom_\2C(\gamma_k\rho,\sigma) \,
  \bigotimes \, \Hom_\2E(\iota,\gamma_k). \label{arr3}\ee
Thus the subspaces $\Hom_\2E(\rho,\sigma)\subset\Hom_{\2C\rtimes\2S}(\rho,\sigma)$ for all 
$\rho,\sigma\in\2C$ are determined by the subspaces 
$\Hom_\2E(\iota,\gamma_k)\subset\Hom_{\2C\rtimes\2S}(\iota,\gamma_k)$. \label{crucial2}\elemma
\prf In the entire proof let $\rho,\sigma\in\2C$ be fixed. With the above identification
 of $\2H_k$ and $\Hom_{\2C\rtimes\2S}(\iota,\gamma_k)$ we can rewrite (\ref{arr}) as
\be \Hom_{\2C\rtimes\2S}(\rho,\sigma)=\bigoplus_{k\in\hat{G}} \
   \Hom_\2C(\gamma_k\rho,\sigma) \ \bigotimes\ \Hom_{\2C\rtimes\2S}(\iota,\gamma_k).  \label{arr1}\ee
If $\id_{\gamma_k}\otimes\psi\in\Hom_{\2C\rtimes\2S}(\iota,\gamma_k)$ is contained in 
$\Hom_\2E(\iota,\gamma_k)$ and $S\in\Hom(\gamma_k\rho,\sigma)$ then 
$S\otimes\psi=S\otimes\Omega\mcirc\id_{\gamma_k}\otimes\psi\in\Hom_\2E(\rho,\sigma)$, since 
$S\in\Hom\,\2C\subset\Hom\,\2E$. Thus in (\ref{arr3}) we have the inclusion $\supset$. 
Now we define positive definite scalar products $\langle \cdot,\cdot\rangle_k$ on 
$\Hom_\2C(\gamma_k\rho,\sigma)$ for all $k\in\hat{G}$ as follows:
\be S, T \ \mapsto \
   \langle S,T\rangle_k \,\id_{\gamma_k}= \id_{\gamma_k}\times\ol{R}_\rho^* \mcirc
   (S^*\circ T)\times\id_{\ol{\rho}} \mcirc\id_{\gamma_k}\times\ol{R}_\rho. \ee
We have used that $\gamma_k$ is irreducible, thus 
$\Hom(\gamma_k,\gamma_k)\cong\7C\,\id_{\gamma_k}$.
(Positive definiteness is seen as follows: $\langle S,S\rangle=0$ implies 
$\id_{\gamma_k}\times\ol{R}_\rho^* \mcirc (S^*\circ S)\times\id_{\ol{\rho}} \mcirc\id_{\gamma_k}\times\ol{R}_\rho=0$. By positivity of the $*$-operation this implies
$S\times\id_{\ol{\rho}}\mcirc\id_{\gamma_k}\times\ol{R}_\rho=0$ and using the conjugate
equation this entails $S=0$.) For every $k\in\hat{G}$ pick an orthonormal basis
$\{ W^k_i,\ i=1,\ldots,\dim\Hom(\gamma_k\rho,\sigma)$ in the Hilbert spaces 
$\Hom(\gamma_k\rho,\sigma)$. Every $\tilde{S}\in\Hom_{\2C\rtimes\2S}(\rho,\sigma)$ is of the 
form $\tilde{S}=\oplus_{l\in\hat{G}} \sum_j S^l_j\otimes\psi^l_j$, where 
$S^l_j\in\Hom(\gamma_l\rho,\sigma)$ and $\psi^l_j\in\2H_l$.
Using the above discussion this can be expressed as
\bea \tilde{S} &=& \bigoplus_{l\in\hat{G}} \sum_j \left( \sum_i \langle W^l_i, 
  S^l_j \rangle \,W^l_i \right) \ \bigotimes \ \psi^l_j  \nn\\ 
   &=& \sum_{l\in\hat{G}} \sum_j \sum_i \langle W^l_i, S^l_j \rangle\  W^l_i \mcirc
   \id_{\gamma_l}\otimes\psi^l_j \,\times\,\id_\rho \nn\\ 
   &=& \sum_{l\in\hat{G}} \sum_j \sum_i W^l_i \mcirc \left( \id_{\gamma_l}\,\times\,
  \ol{R}^*_{\rho}\mcirc(W_i^{l^*}\circ S^l_j) \,\times\,\id_{\ol{\rho}}\mcirc 
   \id_{\gamma_l}\times\ol{R}_\rho \right)\mcirc\id_{\gamma_l}\otimes\psi^l_j\,\times
   \, \id_\rho  \nn\\
   &=& \sum_{k,l\in\hat{G}} \sum_j \sum_i  W^l_i \mcirc \left( \id_{\gamma_l}\times
  \ol{R}^*_{\rho} \mcirc (W_i^{l^*}\circ S^k_j) \times\id_{\ol{\rho}}\mcirc 
   \id_{\gamma_k}\times\ol{R}_\rho \right)\mcirc \id_{\gamma_k}\otimes\psi^k_j \times
   \id_\rho  \nn\\
  &=& \sum_{l\in\hat{G}} \sum_i W^l_i\mcirc \left( \id_{\gamma_l}\times\ol{R}^*_\rho 
   \mcirc (W_i^{l*} \circ \tilde{S})\times\id_\rho \right)  \nn\\
  &=& \bigoplus_{l\in\hat{G}} \sum_i W^l_i \bigotimes \Psi^l_i, 
\eea
where 
\be \Psi^l_i=\id_{\gamma_l}\times\ol{R}^*_\rho\ \circ\ ({W_i^l}^*\circ\tilde{S})\times
   \id_\rho \ \circ \ \ol{R}_\rho\ \in\Hom_{\2C\rtimes\2S}(\iota,\gamma_k). \ee
In the second step we have used 
$S\otimes\psi^l=S\mcirc \psi^l\times\id_\rho$. The fourth equality
is true since the big bracket is in $\Hom_\2C(\gamma_k,\gamma_l)$, which vanishes for
$k\ne l$. In the fifth step we used the interchange law (in $\2C\rtimes\2S$) as in
the following diagram and performed the summations over $k$ and $j$. Now we have
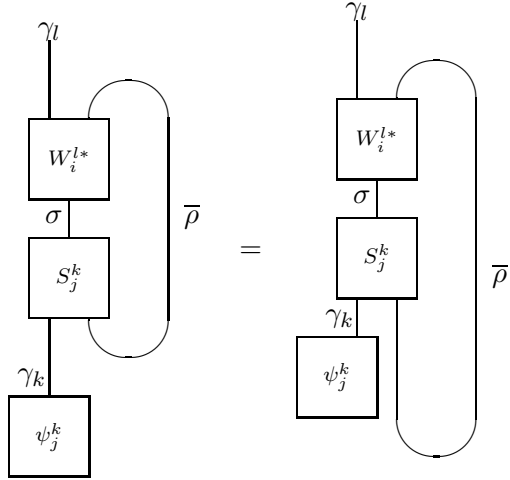
\[ \ba{ccc}  
\begin{tangle} 
\step\object{\gamma_l}\\
\step\id\step\coev \\
\step\Frabox{W_i^{l*}}\Step\id \\
\hh\hstep\step[.4]\mobj{\sigma}\step[.6]\id\step[2.5]\id\step[.4]\mobj{\ol{\rho}} \\
\step\Frabox{S_j^k}\Step\id \\
\step[.2]\mobj{\gamma_k}\step[.8]\id\step\ev\\
\hstep\Frabox{\psi^k_j}
\end{tangle} 
 & \enspace\enspace = \enspace\enspace &
\begin{tangle} 
\step\object{\gamma_l}\\
\step\id\step\coev \\
\step\Frabox{W_i^{l*}}\Step\id \\
\hh\hstep\step[.4]\mobj{\sigma}\step[.6]\id\step[2.5]\id \\
\step\Frabox{S_j^k}\Step\id\step[.4]\mobj{\ol{\rho}} \\
\hh\step[.2]\mobj{\gamma_k}\step[.8]\id\step\id\Step\id \\
\Frabox{\psi^k_j}\step\id\Step\id \\
\Step\ev
\end{tangle} 
\ea \]
and if $\tilde{S}\in\Hom_\2E(\iota,\gamma_l)$ then also
$\Psi^l_i\in\Hom_\2E(\iota,\gamma_l)$ since $W^l_i$ and $\ol{R}_\rho$ are morphisms in $\2C$,
thus in $\2E$. This proves the inclusion $\subset$ in (\ref{arr3}). \qed

\bprop Let $\2E$ be a T$C^*$ such that $\2C\subset\2E\subset\2C\rtimes\2S$. Then 
$\2E=(\2C\rtimes\2S)^H$ where $H=\aut_\2E(\2C\rtimes\2S)$ is a closed subgroup of 
$G=\aut_\2C(\2C\rtimes\2S)$. \label{P1}\eprop
\prf Let $\2F$ be the full subcategory of $\2E$ defined by 
$\obj\,\2F=\obj\,\2E\cap\obj\,\2S\rtimes\2S$, i.e.\ $(\rho,E)\in\2E$ is in $\2F$ iff
$\rho\in\2S$. Then we have the following diagram:
\be\ba{ccccc}
  \2C & \subset & \2E & \subset & \2C\rtimes\2S \\
  \cup &  & \cup &  & \cup \\
  \2S & \subset & \2F & \subset & \2S\rtimes\2S 
\ea\ee
Here all vertical inclusions are full and all categories in the lower row are symmetric.
($\2S\rtimes\2S$ is symmetric since it is the closure under subobjects of 
$\2S\rtimes_0\2S$. The latter is a symmetric tensor category since $\2S$ -- though not 
necessarily contained in $\2D(\2C)$ -- is trivially contained in $\2D(\2S)=\2S$,
entailing that the symmetric braiding of $\2S$ lifts to $\2S\rtimes_0\2S$.)
Fixing a DR representation functor $F: \2S\rightarrow\2H$, where $\2H$ is the 
symmetric tensor category of finite dimensional Hilbert spaces, we define $G$ to be the
group of natural automorphisms of $F$ and have
$\aut_\2S(\2S\rtimes\2S)\cong\aut_\2C(\2C\rtimes\2S)\cong G$. 
Defining $H=\aut_\2F(\2S\rtimes\2S)\subset G$, the proposition follows easily as soon as
we prove
\be \2F=(\2S\rtimes\2S)^H \label{xyz}\ee
since this implies $\Hom_\2F(\iota,\gamma)=\Hom_{\2S\rtimes\2S}(\iota,\gamma)^H, \gamma\in\2S$ 
and by Lemma \ref{crucial2} we have 
$\Hom_\2E(\rho,\sigma)=\Hom_{\2C\rtimes\2S}(\rho,\sigma)^H$ for all $\rho,\sigma\in\2C$.
Since $\2E$ is supposed closed under subobjects this implies $(\rho,E)\in\obj\,\2E$
if the projection $E\in\Hom_{\2C\rtimes\2S}(\rho,\rho)$ is $H$-invariant. On the other hand, 
$(\rho,E)\in\obj\,\2E$ implies $E\in\Hom\,\2E$ since $E=\id_{(\rho,E)}$ and $\2E$
is a category. Thus $(\rho,E)\in\obj\,\2E$ iff $E\in\Hom_{\2C\rtimes\2S}(\rho,\rho)^H$ and 
therefore $\2E=(\2C\rtimes\2S)^H$. Thus we are left with the proof of (\ref{xyz}). 

Choose a section $\{\gamma_k, k\in\hat{G}\}$ of irreducibles in $\2S\cong U(G)$. 
We begin by showing that $F$ extends to a functor
$\hat{F}: \2S\rtimes_0\2S\rightarrow\2H$. For $S\in\Hom(\gamma_k\rho,\sigma), \psi\in\2H_k$
we recall that $S\otimes\psi\in\Hom_{\2S\rtimes_0\2S}(\rho,\sigma)$ and define
$\hat{F}(S\otimes\psi): F(\rho)\rightarrow F(\sigma)$ by
\be \hat{F}(S\otimes\psi)(\phi)=F(S)(\psi\boxtimes\phi), \quad \phi\in F(\rho). \ee
This makes sense since $\psi\boxtimes\phi\in F(\gamma_k)\boxtimes F(\rho)$ and
the latter Hilbert space is canonically isomorphic to $F(\gamma_k\rho)$. By definition
$\hat{F}$ coincides with $F$ on the objects, and it is easy to see that the same is true
on the morphisms $\Hom_\2S(\rho,\sigma)$ of $\2S$.
We have to show that $\hat{F}$ is a symmetric tensor $*$-functor, i.e.\ compatible with
the operations $\circ,\times, *$. We do this only for $\circ$ and leave the other 
arguments to the reader. Let $S\in\Hom(\gamma_k\sigma,\eta), T\in\Hom(\gamma_l\rho,\sigma),
\psi_k\in\2H_k, \psi_l\in\2H_l$ and $\phi\in\2H_\rho=F(\rho)$. We have to show that 
\be \hat{F}(S\otimes\psi_k\mcirc T\otimes\psi_l)\,\phi=
   \hat{F}(S\otimes\psi_k)\circ\hat{F}(T\otimes\psi_l)\,\phi \quad\forall\phi\in\2H_\rho.
\ee
The right hand side equals 
\be \hat{F}(S\otimes\psi_k)F(T)(\psi_l\boxtimes\phi)=
  F(S)(\psi_k\boxtimes F(T)(\psi_l\boxtimes\phi))=F(S\circ\id_{\gamma_k}\times T)
   (\psi_k\boxtimes\psi_l\boxtimes\phi), \ee
and is seen to coincide with the left hand side 
\bea \lefteqn{\hat{F}\left( \bigoplus_{m\in\hat{G}} \sum_{\alpha=1}^{N_{k,l}^m}
   S\,\circ \,\id_{\gamma_k}\times T \,\circ\,
   V_{k,l}^{m,\alpha}\times\id_\rho \ \bigotimes \ F(V_{k,l}^{m,\alpha})^* 
   (\psi_k\boxtimes\psi_l) \right) \phi=}  \\
  &&= \bigoplus_{m\in\hat{G}} \sum_{\alpha=1}^{N_{k,l}^m} F(S\,\circ \,\id_{\gamma_k}
   \times T \,\circ\, V_{k,l}^{m,\alpha}\times\id_\rho) (F(V_{k,l}^{m,\alpha})^* 
   (\psi_k\boxtimes\psi_l)\boxtimes\phi) \nn\eea
appealing to the completeness relation for the bases $\{ V_{k,l}^{m,\alpha} \}$.
The extension of $\hat{F}$ to the new objects $(\rho,E), E\lneq\id_\rho$ 
of $\2S\rtimes\2S=\ol{\2S\rtimes_0\2S}$ is obvious:
$\hat{F}((\rho,E))=\hat{F}(E)\hat{F}(\rho)$, the right hand side being a subspace of the
Hilbert space $\hat{F}(\rho)=F(\rho)$. The functor $\hat{F}: \2S\rtimes\2S\rightarrow\2H$
thus obtained is a symmetric tensor $*$-functor and thus a DR representation functor. 
Furthermore, $\hat{F}\restr\2F$ is a representation functor for $\2F$, and
$\gal(\2F)$ is the set of natural transformations of $\hat{F}\restr\2F$,
i.e.\ the set of families of unitary maps 
$\{ U_{(\rho,E)}\in F(\Hom((\rho,E),(\rho,E))),\ (\rho,E)\in\2F \}$ such that 
\be U_{(\sigma,F)} \circ F(\tilde{S})=F(\tilde{S}) \circ U_{(\rho,E)} \label{comm}\ee
for all $(\rho,E),(\sigma,F)\in\2F, \tilde{S}\in\Hom_{\2F}((\rho,E),(\sigma,F))$. 
Since $\2F$ contains $\2S$, a natural transformation of $\hat{F}\restr\2F$ 
restricts to one of $F$: $\{ U_{(\rho,\id_\rho)},\ \rho\in\2S \}$. 
Now, the group of natural automorphisms of $F$ is just the Galois group $G=\gal(\2C)$.
Let $g\in G$ and let $\{U_{(\rho,\id_\rho)}=\pi_\rho(g), \rho\in\2S\}$ be the 
corresponding natural transformation. A necessary condition for the latter to arise from
a natural transformation of $\hat{F}\restr\2F$ is that (\ref{comm}) holds for all
$\rho, \sigma, \tilde{S}\in\Hom_{\2F}(\rho,\sigma)$. The corresponding $g\in G$ clearly 
constitute a subgroup $H\subset G$. In order to study this subgroup let 
$\tilde{S}\in\Hom_{\2F}(\rho,\sigma)\subset\Hom_{\2S\rtimes\2S}(\rho,\sigma)$. With
\be \tilde{S}=\bigoplus_{k\in\hat{G}} \sum_i \ S_k^i\bigotimes \psi_k^i,\quad
   S_k^i\in\Hom(\gamma_k\rho,\sigma),\ \psi_k\in\2H_k \ee
and the definition of $\hat{F}$ we have 
\be \hat{F}\left(\bigoplus_{k\in\hat{G}} \sum_i \ S_k^i\bigotimes \psi_k^i
  \right)\phi= \sum_{k\in\hat{G}}\sum_i  F(S_k^i)(\psi_k^i\boxtimes\phi). \ee
Then (\ref{comm}) takes the form
\bea \sum_{k\in\hat{G}}\sum_i F(S_k^i)(\psi_k^i\boxtimes\pi_\rho(g)\phi) 
  &=& \pi_\sigma(g)\sum_{k\in\hat{G}}\sum_i F(S_k^i)(\psi_k^i\boxtimes\phi) \nn\\
  &=& \sum_{k\in\hat{G}}\sum_i F(S_k^i)(\pi_k(g)\psi_k^i\boxtimes\pi_\rho(g)\phi). \eea
Since the subspaces $\Hom(\gamma_k\rho,\sigma) \bigotimes \2H_k$ for different $k$ are 
linearly independent, this is true iff $\alpha_g(\tilde{S})=\tilde{S}$. Since this 
must hold for all arrows $\tilde{S}$ in $\2F$ we define
\be H=\{ g\in G\ |\ \alpha_g(\tilde{S})=\tilde{S} \quad \forall \tilde{S}\in \2F\}, 
\ee
which is a closed subgroup of $G$.
For $g\in H,\ U_{(\rho,\id_\rho)}=\pi_\rho(g)$ commutes with the projections
$E\in\Hom_{\2F}(\rho,\rho)$, and $U_{(\rho,E)}=U_{(\rho,\id_\rho)}\restr E\2H_\rho$ is a
natural transformation of $\hat{F}\restr\2F$. Thus $\gal(\2F)\cong H$, and by
the duality theorem we know that $\2F$ is a category of representations of $H$.
Thus for $T\in\Hom_{\2S\rtimes\2S}(\rho,\sigma)$ the linear operator 
$\hat{F}(T): F(\rho)\rightarrow F(\sigma)$ is contained in 
$\hat{F}(\Hom_{\2F}(\rho,\sigma))$ iff it intertwines the representations $\pi_\rho$ and
$\pi_\sigma$. By the above this is equivalent to $T$ being $H$-invariant and therefore
we have $\Hom_\2F(\rho,\sigma)=\Hom_{\2S\rtimes\2S}(\rho,\sigma)^H$ for $\rho,\sigma\in\2S$.
For the subobjects $(\rho,E)$ the argument at the beginning of the proof applies and we
obtain $\2F=(\2S\rtimes\2S)^H$. \qed

Now we consider the question for which subgroups $H\subset G$ there is a subcategory 
$\2T\subset\2S$ such that $(\2C\rtimes\2S)^H\cong\2C\rtimes\2T$. We begin with 
three lemmas.

\blemma Let $G$ be a compact group and let $\pi$ be an irreducible unitary 
representation on the Hilbert space $\2H$. Let $H$ be a closed normal subgroup of
$G$. Then the subspace $\2H^H\subset\2H$ of $H$-invariant vectors is either $\{0\}$ or
$\2H$. \label{L1}\elemma
\prf Let $\psi\in\2H$ be $H$-invariant. The normality of $H$ implies that the vectors
$\pi(g)\psi,\ g\in G$ are $H$-invariant, too. But the span of the latter is $\2H$, since
otherwise it would be a non-trivial $G$-invariant subspace, which does not exist by 
irreducibility of $\pi$. \qed

\blemma Let $G$ be compact and $H$ be a closed normal subgroup. Then there is a 
one-to-one correspondence between the (i) continuous unitary representations $\pi$ of 
$G/H$ and (ii) continuous unitary representations $\hat{\pi}$ of $G$ such that 
$H\subset\ker\,\hat{\pi}$. This correspondence restricts to irreducible representations.
An intertwiner between representations $\pi, \pi'$ lifts to $\hat{\pi}, \hat{\pi}'$
and vice versa. \label{L2}\elemma
\prf Let $\phi: G\rightarrow G/H$ be the quotient homomorphism. Then the correspondences
are given by $\pi\mapsto\hat{\pi}=\pi\circ\phi$ and 
$\hat{\pi}\mapsto\pi=\hat{\pi}\circ\phi^{-1}$, where the latter is well-defined since
$\hat{\pi}$ is constant on cosets. These constructions respect continuity since $\phi$
is continuous and open. The statement on intertwiners is obvious. \qed

The following is not explicitly contained in \cite{dr1}, but a part of the results is
contained in the more general \cite[Thm.\ 6.10]{dr1}. 
\blemma Let $\2S$ be a ST$C^*$ with $\gal(\2S)\cong G$. Pick a representation 
functor $F$ of Doplicher and Roberts which identifies $\2S$ with a category $U(G)$ of 
representations of $G$ and let $\pi_\rho$ be the action of $G$ on the Hilbert space
$F(\rho)$. For a closed normal subgroup $H$ of $G$ the full subcategory of $\2S$ 
defined by $\obj\,\2T_H=\{\rho\in\2S \ | \ H\subset\ker\,\pi_\rho\}$ is a 
replete full symmetric subcategory with conjugates etc., and $\gal(\2T_H)\cong G/H$. 
The map $H\mapsto\2T_H$ is bijective, the inverse being given by
$\2T\mapsto H_\2T=\{h\in G\ |\ h\in\ker\,\pi_\rho\ \forall \rho\in\2T\}$. 
(In these considerations the non-uniqueness of the
functor $F$ is unimportant since the kernel of the representation $\pi_\rho$
does not depend on the choice of $F$.) \label{L3}\elemma
\prf Given a closed normal subgroup $H$, define $\2T_H\subset\2S$ as given. The braiding
and the $*$-operation restrict to $\2T_H$, which is also closed under conjugates, direct 
sums and subobjects. For $\rho\in\2T_H$ Lemma \ref{L2} gives rise to a representation of 
$G/H$ on $F(\rho)$, and $F(T)$ where $\rho,\sigma\in\2T_H, T\in\Hom(\rho,\sigma)$ intertwines
the representations of $G/H$ on $F(\rho), F(\sigma)$. Since $U(G)=F(\2S)$ is complete in
the sense that for every $g\in G$ there is a $\rho\in\2S$ such that
$F(\rho)(g)\ne\11$, the same holds for $\2T_H$ and $G/H$, which implies
$\gal(\2T_H)\cong G/H$.
On the other hand, given $\2T,\ H_\2T$ clearly is a closed normal subgroup of $G$, and 
we have to show that this map is inverse to $H\mapsto \2T_H$. Obviously, 
$H\subset H_{\2T_H}$ 
and $\2T\subset \2T_{H_\2T}$. By the above, $F(\2T_H)$ can be looked at as a complete 
category of representations of $G/H$. Thus $g\in G$ is in $H_{\2T_H}$ iff $\ol{g}=\ol{e}$
(where $\ol{g}=\phi(g)$ is the image of $g$ in $G/H$) iff $g\in H$, whence
$H_{\2T_H}= H$. For given $\2T\subset\2S$, $F$ restricts to an embedding functor for 
$\2T$, and $\gal(\2T)$ is (isomorphic to) the group of natural transformations of
$F\restr\2T$ to itself. Since $g\in\gal(\2S)$ is trivial as a natural transformation of
$F\restr\2T$ iff $g\in H_\2T$ we have a homomorphism of $G/H_\2T$ into $\gal(\2T)$.
Since the map $\phi:G\rightarrow G/H_\2T$ is surjective we have in fact an isomorphism
$\gal(\2T)\cong G/H_\2T$. Comparing this with
$\gal(\2T_{H_\2T})\cong G/H_{\2T_{H_\2T}}=G/H_\2T$, where we have used $H_{\2T_H}=H$,
this implies $\2T\simeq\2T_{H_\2T}$. Since $\2T, \2T_{H_\2T}$ are replete full
subcategories of $\2S$ we have $\2T=\2T_{H_\2T}$. \qed

\bprop Given $\2C\subset\2E\subset\2C\rtimes\2S$ where $\2E\simeq(\2C\rtimes\2S)^H$, the
subgroup $H\subset G$ is normal iff there is a ST$C^*$ $\2T\subset\2S$ such that
$\2E\cong\2C\rtimes\2T$. In this case $\aut_\2E(\2C\rtimes\2S)=H$ and
$\aut_\2C(\2E)\cong G/H$. \label{P2}\eprop
\prf Let $H$ be a normal subgroup of $G$. Pick a functor $F_\2S$ identifying $\2S$ 
with a category $U(G)$ of 
representations of $G$. Let $\2T_H\subset\2S$ be the full subcategory corresponding 
to $H$. $F_\2S$ restricts to $\2T_H$, and when comparing $\2C\rtimes\2S, \2C\rtimes\2T_H$
we will choose the functors $F_\2S, F_\2S\restr\2T_H$ in the construction of the crossed
products.

By definition $(\2C\rtimes\2S)^H\cong\2C\rtimes\2T_H$ is the subcategory of 
$\2C\rtimes\2S$ whose objects and arrows are $H$-invariant. In view of (\ref{arr}) and 
(\ref{alpha1}) this means for $\rho,\sigma\in\2C$ that
\be \Hom_{(\2C\rtimes\2S)^H}(\rho,\sigma) = \bigoplus_{k\in\hat{G}} \
    \Hom(\gamma_k\rho,\sigma) \bigotimes \2H_k^H =
   \bigoplus_{k\in\hat{G} \atop H\subset\ker\pi_k} \
   \Hom(\gamma_k\rho,\sigma) \bigotimes \2H_k,  \label{x1}\ee
where in the second step we have applied Lemma \ref{L1}. On the other hand
\be \Hom_{\2C\rtimes\2T_H}(\rho,\sigma) = \bigoplus_{k\in\widehat{G/H}} \
    \Hom(\gamma_k\rho,\sigma) \bigotimes \2H_k, \label{x2}\ee
where $\2H_k$ now carries an irreducible representation of $G/H$. By Lemma \ref{L2} there
is a canonical one-to-one correspondence between $k\in\widehat{G/H}$ and 
$k\in\hat{G}, H\subset\ker\pi_k$. Choosing the same $\gamma_k's$ in (\ref{x2})
as in (\ref{x1}) we can identify the right hand sides of (\ref{x1}) and (\ref{x2}), and 
the products $\circ, \times$ on the arrows of $(\2C\rtimes\2S)^H$ and $\2C\rtimes\2T_H$
are the same since $F_{\2T_H}$ is the restriction of $F_\2S$ to $\2T_H$. In view of
$\Hom_{(\2C\rtimes\2S)^H}(\rho,\rho)=\Hom_{\2C\rtimes\2T_H}(\rho,\rho)$ also the objects of 
$(\2C\rtimes\2S)^H$ can be identified with those of $\2C\rtimes\2T_H$. Thus 
$(\2C\rtimes\2S)^H\cong\2C\rtimes\2T_H$.
The preceding argument depended on choosing $F_\2S\restr\2T_H$ for the definition of 
$\2D\rtimes\2T_H$. But by Prop.\ \ref{uniq} another choice of $F_{\2T_H}$ yields an 
isomorphic crossed product category.
Conversely, consider $\2T\subset\2S$ where $\gal(\2S)\cong G=\aut_\2C(\2C\rtimes\2S)$.
Then there is  a normal subgroup $H_\2T$ of $G$ such that $\gal(\2T)\cong G/H$, and it 
is easy to verify that $\2C\rtimes\2T\cong(\2C\rtimes\2S)^{H_\2T}$. \qed

The preceding results were of a relative nature, classifying intermediate extensions
$\2E$ such that $\2C\subset\2E\subset\2C\rtimes\2S$, where $\2S\subset\2D$ was not
assumed. The following result clarifies the role of the absolute Galois groups for
extensions $\2C\rtimes\2S$ where $\2S\subset\2D$.
\bprop For $\2S\subset\2D$ we have
$\gal(\2C\rtimes\2S)\cong\aut_{\2C\rtimes\2S}(\2C\rtimes\2D)\cong H_\2S$.
\eprop
\prf By Cor.\ \ref{gcoro} we have $\2D(\2C\rtimes\2S)\cong\2D\rtimes\2S$. That the 
compact group associated to the ST$C^*\ \2D\rtimes\2S$ is $H_\2S$ follows from the proof
of Prop.\ \ref{P1}. \qed\\
\rems 1. In particular, for $\2S=\2D$ we have $\gal(\2S)=\gal(\2C)=\gal(\2D)$ and thus
$\gal(\2C\rtimes\2D)=\11$, which was the statement of Thm.\ \ref{main}. \\
2. Let $\2C$ be symmetric, i.e.\ $\2D=\2C$ with $\gal(\2C)\cong G$.
Then taking the crossed product $\2C\rtimes\2S$ with $\2S\subset\2D$ and 
$\gal(\2S)\cong G/H$ amounts to restricting the representations in $U(G)\cong\2C$ to 
the normal subgroup $H$. Then the statement of Prop.\ \ref{crucial} on the equality of
multiplicities and dimensions is nothing but the well known result \cite[\S 49]{cr}.
Namely given an irreducible representation $\pi$ of $G$ all irreducible representations
of $H$ in $\pi\restr H$ occur with the same multiplicity and have the same dimensions.

The preceding results make the analogy with algebraic field extensions $K\supset F$
obvious. Also these can be iterated until one arrives at the algebraic closure 
$\ol{F}$. The latter is the unique (up to $F$-isomorphism) algebraic extension in which 
all polynomials split into linear factors with the consequence that further algebraic 
extensions do not exist. Furthermore, there is a one-to-one relation between Galois
extensions $K\supset F$ and closed normal subgroups $H$ of the absolute Galois group of
$F$, given by $H\mapsto \ol{F}^H,\ K\mapsto \aut_K\,\ol{F}$.

Observe that the analogy with the algebraic closure is -- of course -- not quite perfect
since $\ol{\ol{\2C}}$ may have less irreducible objects than $\2C$ or be even trivial:
\bprop $\2C\rtimes\2S$ is trivial -- in the sense that all irreducible objects are 
equivalent to the identity object $\iota$ -- iff $\2S=\2D=\2C$. Equivalently, $\2C$ is 
completely degenerate, i.e.\ symmetric, and $\2S=\2C$. \eprop
\prf If $\2S$ is strictly smaller than $\2D$ then $\2C\rtimes\2S$ by the above contains
degenerate objects which are inequivalent to $\iota$. Thus assume $\2S=\2D$.
The irreducible objects of $\2C\rtimes\2D$ are obtained by decomposing those of 
$\2C$. We have seen that the degenerate objects of $\2C$ become multiples of the identity
in $\2C\rtimes\2D$. But the decomposition in $\2C\rtimes\2D$ of a non-degenerate object 
of $\2C$ yields non-degenerate objects, which are inequivalent to $\iota$. \qed
\bcoro $\ol{\ol{\2C}}=\2C\rtimes\2D$ is non-trivial iff $\2C$ contains at least one
non-degenerate object. \ecoro

\sectreset{Further Directions}
\subsection{Abelian Groups $G$}\label{abel}
In this subsection we consider the special case where all irreducible objects in $\2S$ 
have dimension one, which is equivalent to $\gal(\2S)$ being abelian. Our aim will 
be to give an explicit description of the sector structure of $\2C\rtimes\2S$, where a 
{\it sector} is a unitary isomorphism class of objects. (Abusing notation we write
$\gamma, \rho$ etc.\ for objects and for the corresponding sectors.)

Denoting by $\Delta$ the set of irreducible sectors of $\2C$, the tensor product and
braiding in $\2C$ render $\Delta$ an abelian semigroup. $\Delta$ decomposes into the set 
$\Delta_\2S$ of sectors in $\2S$ and the complement $\Delta'$. Under the above 
assumption of one-dimensionality $K\equiv \Delta_\2S$ is a discrete abelian group and 
the compact DR group is just the Pontrjagin dual $G=\hat{K}$. Given an irreducible 
$\gamma\in K$ and an irreducible $\rho\in\Delta'$, 
$\gamma\rho$ is irreducible (in $\2C$) since by $d_\gamma=1$ and Frobenius reciprocity 
$\Hom(\gamma\rho,\gamma\rho)\cong\Hom(\rho,\rho)\cong\7C$. Another use of Frobenius reciprocity 
\cite[Lemma 3.9]{mue4} shows that $\gamma\rho$ is in $\Delta'$. Thus the sectors in $K$
act on those in $\Delta'$ by permutation, which implies that $\Delta'$ decomposes into 
$K$-orbits $\ul{\rho}:= \{ \gamma\rho, \gamma\in K\}$.
Given  irreducible $\rho_1,\rho_2\in\Delta$, (\ref{arr}) implies that 
$\rho_1,\rho_2$ are unitarily equivalent in $\2C\rtimes\2S$ iff
$\rho_1,\rho_2$ are in the same orbit (i.e.\ $\ul{\rho_1}=\ul{\rho_2}$) and disjoint 
otherwise. Thus in order to find all sectors in $\2C\rtimes\2S$ it suffices to consider 
one element $\rho$ of each orbit $\ul{\rho}$ and to decompose it into irreducibles.
Since $\2C\rtimes\2S$ is closed under direct sums and subobjects the decomposition of
$\rho$ is governed by the semi-simple algebra $\Hom_{\2C\rtimes\2S}(\rho,\rho)$. 
It is well-known that 
\be \Hom_{\2C\rtimes\2S}(\rho,\rho) = \bigoplus_{i\in I} M_{N_i} \quad\impl\quad
  \rho=\bigoplus_{i\in I} N_i\,\rho_i. \ee
Here $M_d$ is the full matrix algebra of rank $d$, thus dimension $d^2$, and the $\rho_i$
are pairwise inequivalent irreducible sectors, occurring with multiplicity $N_i$.

Now we work out explicitly the structure of $\Hom_{\2C\rtimes\2S}(\rho,\rho)$.
Motivated by (\ref{compos}) and the fact that the spaces $\Hom(\gamma_k\rho,\rho)$ are either
zero or one dimensional we define
\be K_\rho = \{ k\in K,\ \gamma_k\rho\cong\rho \}, \ee
which clearly is a subgroup of $K$. By remark 1 after Defin.\ \ref{maindef} $K$ 
is finite. Chosing unitary intertwiners 
$T_k\in\Hom(\gamma_k\rho,\rho), k\in K_\rho$ and normalized vectors $\psi_k\in\2H_k$,
$\Hom_{\2C\rtimes\2S}(\rho,\rho)$ is spanned by $\{T_k\otimes\psi_k,\,k\in K_\rho\}$ 
and we have
\be T_k\otimes\psi_k\mcirc T_l\otimes\psi_l = 
    T_k\,\circ \,\id_{\gamma_k}\times T_l \,\circ\,
   V_{k,l}^{kl}\times\id_\rho \ \bigotimes \ F(V_{k,l}^{kl})^*  (\psi_k\boxtimes\psi_l).
\ee
Now
\be T_k\,\circ \,\id_{\gamma_k}\times T_l\,\circ\,  V_{k,l}^{kl}\times\id_\rho
  \in\Hom(\gamma_{kl}\rho,\rho), \ee
and since $\Hom(\gamma_{kl}\rho,\rho)$ is one dimensional we have
$T_k\,\circ\,\id_{\gamma_k}\times T_l\,\circ\, V_{k,l}^{kl}\times\id_\rho\propto T_{kl}$.
Similarly, $F(V_{k,l}^{kl})^*  (\psi_k\boxtimes\psi_l)$ is a unit vector in 
$\2H_{kl}$, thus proportional to $\psi_{kl}$. Therefore
\be T_k\otimes\psi_k\mcirc T_l\otimes\psi_l = c(k,l)\ T_{kl}\otimes\psi_{kl}, 
\label{twga}\ee
where associativity implies $c$ to be a 2-cocycle in $Z^2(K_\rho,\7T)$, and 
$\Hom_{\2C\rtimes\2S}(\rho,\rho)$ is the twisted group algebra $\7C^c K_\rho$.
(This result could also have been derived from the general theory of ergodic actions of
compact abelian groups on von Neumann algebras, cf.\ eg.\ \cite{ahk}.)
Due to $T_e\in\Hom(\rho,\rho)\in\7C\,\id_\rho$ we can choose $T_e=\id_\rho$, which will
always be assumed in the sequel. Now we need some group theoretical results.

\blemma Let $A$ be a finite abelian group and $c\in Z^2(A,\7T)$. Then the center
of the twisted group algebra $\7C^cA=\mbox{span}\{ U_k, k\in A\}$ with 
$U_kU_l=c(k,l)U_{kl}$ is spanned by $\{ U_k,\ k\in B \}$, where
\be B = \{ k\in A \ |\  c(k,l)=c(l,k)\ \forall l\in B \}\ee
is a subgroup of $A$. In fact, $Z(\7C^cA)\cong\7CB\cong\7C(\hat{B})$.
The twisted group algebra $\7C^cA$ is isomorphic to the 
tensor product of its center with a full matrix algebra:
\be \7C^cA \cong M_N \otimes \7C(\hat{B}) \cong
   \underbrace{M_N \oplus M_N \oplus \ldots \oplus M_N}_{|B|\mbox{ terms}}, \ee
where $N=\sqrt{|A|/|B|}$. The minimal projections of the center are labeled by the 
elements of the dual group $\hat{B}$ and under the canonical action of the dual group
$\hat{A}$ they are permuted according to
\be \alpha_g(P_\chi)=P_{\ol{g}\chi}, \ee
where $\ol{g}\in\hat{B}$ is the restriction of the character $g\in\hat{A}$ to the 
subgroup $B\subset A$. \label{untw}\elemma
\prf The twisted group algebra $\7C^cA$ is a von Neumann algebra. This can be shown
by explicitly exhibiting a positive $*$-operation or by considering $\7C^cA$
as a twisted product of the von Neumann algebra $\7C$ with $A$. Since the canonical
action of the dual group $\hat{A}$ on $\7C^cA$ given by 
$\alpha_g(U_k)=\langle g,k\rangle\, U_k$ is ergodic, Lemma \ref{ergod} applies and 
gives the result on the tensor product structure of the twisted group algebra.
The claim on the center follows by specialization to an abelian group $A$ of
well-known results on the center of twisted group algebras, cf.\ \cite{kp}, or by an 
easy direct proof. That $B$ is a subgroup of $A$ is then obvious in view of
(\ref{twga}) and the fact that the center is a subalgebra. Now, in restriction to $B$ 
the cocycle $c$ is symmetric, which is equivalent to $c\restr B$ being a coboundary: 
\be c(k,l)=\frac{f(kl)}{f(k)f(l)} \quad \forall k,l\in B. \ee
With the replacement $U_k\rightarrow f(k) U_k, \, k\in B$ the cocycle 
disappears on $B$ and we have $Z(\7C^cA)\cong\7CB$. By Pontrjagin duality this
is isomorphic to $\7C(\hat{B})$ and the minimal projections in the center 
are given by
\be P_\chi = \frac{1}{|B|}\sum_{k\in B} \chi(k)\,U_k, \label{proj}\ee
where $\chi\in\hat{B}$ is a character of $B$. From this formula it is obvious that the 
action of $\hat{A}$ permutes these projections as stated. \qed

Applying Lemma \ref{untw} to $\rho$ with $A=K_\rho$ 
and $U_k=T_k\otimes\psi_k$, we define $L_\rho$ to be the $B$ of the lemma and obtain
\bprop In $\2C\rtimes\2S$ the object $\rho\in\2C$ decomposes according to
\be \rho\cong N_\rho \, \bigoplus_{\chi\in\widehat{L_\rho}} \rho^\chi, \ee
where the $\rho^\chi,\ \chi\in\widehat{L_\rho}$ are irreducible, mutually 
inequivalent and all occur with the same multiplicity $N_\rho=\sqrt{|K_\rho|/L_\rho|}$.
The automorphism group $G$ of $\2C\rtimes\2S$ permutes the subsectors according to 
$\alpha_g(\rho^\chi)\cong\rho^{\ol{g}\chi}$. Here $\ol{g}\in\widehat{L_\rho}$ is the 
restriction of $g\in G=\hat{K}$, considered as a character on $K$, to the subgroup 
$L_\rho\subset K_\rho\subset K$. \eprop
\rem The result that all irreducible components of $\rho$ appear with the same
multiplicity $N_\rho$ appears as the (unproved) assumption of `fixpoint homogeneity'
in conformal field theory, cf.\ \cite{fss}.

\bcoro The irreducible sectors (isomorphism classes of irreducible objects) of
$\2C\rtimes\2S$ are labeled by pairs $(\ul{\rho},\chi)$. Here 
$\ul{\rho}\in\Delta/\Delta_\2S$ is an orbit of irreducibles in $\Delta$ under the action
of the group $\Delta_\2S$ of degenerate sectors by multiplication and $\chi$ is a 
character of the subgroup $L_{\ul{\rho}}\subset K_{\ul{\rho}}$. \ecoro

\subsection{Remarks on the case $\2S\not\subset\2D$}\label{nobraid}
Whereas the definition of $\2C\rtimes\2D$ does not require $\2S\subset\2D$, we have seen
that only under this condition the braiding $\ve$ of $\2C$ gives rise to a braiding for
$\2C\rtimes\2S$. Even though this was without importance for the larger part of Sect.\ 
\ref{sect4} we remark that also in the case $\2S\not\subset\2D$ one can obtain braided
tensor categories, which is of relevance for the applications to conformal quantum field
theory, in particular the theory of modular invariants, as well as to subfactor theory.

If $\2S\not\subset\2D$ we can still obtain a {\it braided} semidirect product if we 
replace $\2C$ by the replete full subcategory $\2C_\2S$ which is defined by
\be \obj\,\2C_\2S=\{\rho\in\2C\ | \ \ve_M(\rho,\gamma)=\id_{\rho\gamma}\quad 
  \forall \gamma\in\2S \}. \ee
This set is easily seen to be closed under isomorphism, tensor products, conjugates,
direct sums and subobjects. Since $\2S$ is symmetric we clearly have $\2C_\2S\supset\2S$,
and by definition $\2S\subset\2D(\2C_\2S)$. 
Thus $\2C_\2S$ satisfies all assumptions and we can construct $\2C_\2S\rtimes\2S$, 
which is a non-trivial braided tensor category unless $\2C_\2S=\2S$. 
(It may be instructive to compare $\2D(\2C)$ with the center $Z(M)$ of a von 
Neumann algebra $M$, $\2S$ with an abelian subalgebra $A\subset M$ and $\2C_\2S$ with the
relative commutant $M\cap A'$. Then $\2C_\2S=\2S$ corresponds to $M\cap A'=A$, i.e.\
$A$ maximal abelian in $M$.)
By the preceding discussion $\2C_\2S\rtimes\2S$ will be non-degenerate iff 
$\2S=\2D(\2C_\2S)$, which can of course be enforced by replacing $\2S$ by $\2D(\2C_\2S)$.
This makes clear that given a pair $(\2C,\2S)$ where $\2C$ is a BT$C^*$ with $\2S$ a 
symmetric subcategory and setting
\be \2C'= \2C_\2S,\quad \2S'= \2D(\2C_\2S) \ee
we obtain a non-degenerate BT$C^*$ $\2C'\rtimes\2S'$. It would be interesting to
understand the structure of the set of all such crossed products obtainable from a 
given $\2C$.

\subsection{The Case of Supergroups}\label{super}
Up to now we have assumed that all objects in $\2S$ are bosons, i.e.\ have twist equal
to $+1$. Now we consider the general case, assuming that there is at least one fermionic 
degenerate sector. Clearly we may apply the construction as expounded so far to replace 
$\2C$ by $\2C'=\2C\rtimes\2D_+$, where $\2D_+\subset\2D$ is the category of bosonic
degenerate objects. By the above it is clear 
that $\gal(\2C\rtimes\2D_+)\cong\7Z_2$, i.e.\ this BT$C^*$ has only one degenerate 
sector $\gamma$, which satisfies $\gamma^2\cong\iota$ and 
$\ve(\gamma,\gamma)=-\id_{\gamma^2}$.

\blemma A fermionic degenerate object $\gamma$ of dimension one does not have fixpoints,
i.e.\ there is no irreducible $\rho\in\2C$ such that $\gamma\rho\cong\rho$. \elemma
\prf Assume $\rho$ is irreducible such that $\gamma\rho\cong\rho$. Then 
$\omega(\rho)=\omega(\gamma\rho)$. On the other hand, in view of 
$\ve_M(\gamma,\rho)=\id_{\gamma\rho}$ (\ref{tw2}) implies 
$\omega(\rho)=\omega(\gamma\rho)\omega(\gamma)$. Since $|\omega(\rho)|=1$ this is
possible only if $\omega(\gamma)=1$. \qed

Thus $\obj\,\2C'$ decomposes into orbits of length two under the action of $\gamma$ by
multiplication. Assuming naively that as in the bosonic case there is a similar cross 
product construction, which we call $\2C'\rtimes\gamma$, we expect that the 
irreducible objects in $\2C'$ remain irreducible in $\2C'\rtimes\gamma$.
The only effect of the cross product construction should be pairwise identifying the 
objects $\rho$ and $\gamma\rho$ for all $\rho$. The question is whether 
$\2C'\rtimes\gamma$ exists as a BT$C^*$. Unfortunately, this is impossible, since 
$\rho$ and $\gamma\rho$ are equivalent in the would-be BT$C^*\ \2C'\rtimes\gamma$, 
but they have different twist.

This does, of course, not exclude the possibility that there is a full subcategory which 
contains precisely one object from each orbit $\{ \rho, \gamma\rho \}$. But we do not
have a criterion which would guarantee this.

\sectreset{Conclusions and Outlook}
If symmetric tensor categories are considered as an extreme species of braided tensor
categories then non-degenerate categories are the opposite extreme and the construction
of the modular closure $\ol{\ol{\2C}}$ amounts to dividing out the symmetric part. Thus 
$\ol{\ol{\2C}}$ should be considered as the 1-dimensional analogue (in the sense of 
higher category theory) of the quotient group $G/Z(G)$, which for nice $G$ (e.g.\ semisimple)
has trivial center. The significance of our results lies in showing that 
every braided T$C^*$ (braided tensor category plus some additional structure) can be 
faithfully embedded into a braided tensor category which is non-degenerate. In this way 
we obtain a non-trivial category whenever the original one is not symmetric. In 
particular we obtain a unitary (in the sense of \cite{t}) modular category whenever 
$\ol{\ol{\2C}}=\2C\rtimes\2D$ is rational. Since modular categories are instrumental in
the construction of 3-manifold invariants \cite{t} our construction has obvious 
applications to topology.  

Our strategy for removing the degeneracy was to add morphisms to the category $\2C$ and
to close the category $\2C\rtimes_0\2S$ thus obtained w.r.t.\ subobjects. This is 
precisely the approach conjectured to work in \cite[p.\ 460]{tw}:
\begin{quote} 
... it seems likely that one could get more modular categories by adding additional
morphisms to the categories which are constructed in this paper...
\end{quote}
These authors did not, however, indicate a general procedure. Comparison of our
construction in Sect.\ \ref{sect3} with \cite{mue4} reveals that we have done little more
than to translate the formulae {\it derived} in the QFT framework \cite{khr2,mue4} into 
a more abstract setting and re-prove facts like associativity which are obvious in the 
QFT case. The considerations of Sect.\ \ref{sect4}, however, have little in common with 
those \cite{mue4} in the QFT setting. 

Concerning the special case of Subsect.\ \ref{abel} where all irreducible degenerate 
objects have dimension one we cite \cite[p.\ 359]{ker}:
\begin{quote} 
... In case $\2J_0$ is a subgroup of invertibles $\{\sigma\}$ we have for the
natural action of its elements on $k^\2J$ that $\2S\sigma=\2S$. Hence $\2S$ and $\2T$
can be defined on the orbit space $im(\sum_{\sigma\in\2J_0} \sigma)$, where we can hope
for the modularity condition to hold.
\end{quote}
Also this conjecture has been proved above, but as we have seen the decomposition into 
irreducibles of the objects in $\2C\rtimes\2S$ is not quite trivial, since it may be 
complicated by (i) the existence of the stabilizers $K_\rho$ and by (ii) non-trivial 
2-cocycles which lead to multiplicities $N_\rho>1$, cf.\ also \cite{fss}.

We have formulated our results in terms of $C^*$-tensor categories since they are the
natural language for investigations on operator algebras and \qfts. But it should be
clear -- as already pointed out -- that the $C^*$-structure does not play a crucial
role. Replacing the DR duality theorem by the one of Deligne \cite{del2} one can
formulate versions of our construction for braided tensor categories which are enriched
over ${\bf\mbox{Vect}}_k$ for an arbitrary field $k$ of characteristic zero. Also the 
strictness of the tensor categories assumed in this paper is not crucial.
But note that Deligne has to assume integrality of dimensions in the symmetric category,
whereas in the framework of $C^*$-categories this is automatic \cite[Cor.\ 2.15]{dr1}
as a consequence of positivity.

We close by listing some questions which were not treated in this work and directions
for further investigations:
\begin{enumerate}
\item Find a universal property which characterizes the modular closure $\ol{\ol{\2C}}$
up to equivalence.
\item Let $\2C$ be a BT$C^*$, acted upon by a compact group $G$. Under which 
condition is there a subcategory $\2S\subset\2D(C^G)$ with $\gal(\2S)\cong G$ such that
$\2C\simeq\2C^G\rtimes\2S$?
\item Clarify the decomposition of an irreducible $\rho\in\2C$ into subobjects in
$\2C\rtimes\2S$, extending the considerations in Subsect.\ \ref{abel} to the non-abelian
case. \\
This looks difficult since not even the results of Subsect.\ \ref{abel} for the abelian 
case are very explicit. We indicate a generalization of the considerations given above
which, however, is not quite sufficient. Assume $\rho\in\2C$ irreducible is such that 
$\gamma\rho\cong d_\gamma\rho$ whenever $\gamma\in\2S$ is irreducible and 
$\Hom(\gamma\rho,\rho)\ne\{0\}$. (This clearly includes the case of $G$ being abelian.)
Then the set $\{ \gamma\in\2S\ |\ \gamma\rho\cong d_\gamma\rho \}$ is closed under 
multiplication and gives rise to a full subcategory $\2S_\rho$ of the ST$C^*\ \2S$,
which is the representation category of a quotient $G_\rho$ of $G$. Clearly, the action
of $G$ on $\Hom_{\2C\rtimes\2S}(\rho,\rho)$ factors through $G_\rho$. This action of $G_\rho$
has full multiplicity in the sense that the spectral subspace corresponding to any 
irreducible representation $\pi_k$ has dimension $d_k^2$. Then the considerations of
\cite[II]{wa} apply and we know that 
$\Hom_{\2C\rtimes\2S}(\rho,\rho)\cong\pi_\omega(\widehat{G_\rho})$ where $\omega$ is a
2-cocycle on $\widehat{G_\rho}$ and the isomorphism intertwines the actions of $G_\rho$. 
See \cite[II]{wa} for the terminology. This case seems, however, too special 
to deserve further analysis.
\item Given a rational BT$C^*\ \2C$ find a direct construction of the 3-manifold 
invariant arising from the modular closure $\ol{\ol{\2C}}$, bypassing the construction 
of the latter.
\item There is an obvious connection between the crossed product $\2C\rtimes\2S$ and the 
`orbifold constructions in subfactors' \cite{ek,x} which deserves to be worked out.
\item Generalize everything in this paper to the non-connected case where 
$\Hom(\iota,\iota)\ne\7C\,\id_\iota$ and the compact (super)groups are replaced by 
compact (super)groupoids \cite{baez}. The resulting Galois theory should resemble the 
Galois theory for commutative rings instead of the one for fields.
\item Since Janelidze's general Galois theory for categories \cite{j} was modeled on the
Galois theory for commutative rings as expounded by Magid, it should be possible to show 
that with the proper identifications our Galois correspondence fits into Janelidze's formalism, 
also after extension to the non-connected case.
\end{enumerate}

\vspace{2cm}
\noindent{\it Acknowledgments.} This paper would not have been written without 
K.-H.\ Rehren's works \cite{khr1,khr2} and his proposal to pursue the line of research
which led to \cite{mue4}. I am grateful to J.\ Roberts and R.\ Longo for many useful 
discussions and to L.\ Tuset for a discussion on Prop.\ \ref{twist}.
I thank D.\ E.\ Evans for the opportunity to present this work in the Wales
Video Seminar at an early stage and R.\ Brown, T.\ Porter for interest and suggestions 
provided on this occasion. 
The tangle diagrams were generated using the macro package {\tt tangle-s.sty} written by 
Yu.\ Bespalov and V.\ Lyubashenko.

\vspace{1.5cm}
\noindent{\it Note added.} 
After this paper was completed I received the preprint {\it Cat\'{e}gories 
pr\'{e}modul\-aires, modularisations et invariants de vari\'{e}t\'{e}s de dimension 3} by
A.\ Brugui\`{e}res, which was finished several months earlier. In this paper a
construction is given which is equivalent to our definition of $\2C\rtimes\2S$ if
$\2S$ is rational, i.e.\ $\gal(\2S)$ is finite, and $\2S\subset\2D$. Brugui\`{e}res 
does not consider the cases where $\2S\not\subset\2D$ or $|\gal(\2S)|=\infty$, nor does 
he obtain the results of Prop.\ \ref{crucial} and the Galois correspondence. On the 
other hand his construction is more elegant and canonical -- yet less elementary -- in 
that it neither uses a section $\{\gamma_k,\ k\in\hat{G}\}$ nor bases in the intertwiner
spaces, and he solves the problems 1 and 4 listed above. Brugui\`{e}res' work relies on 
Deligne's characterization \cite{del2} of representation categories, which confirms our
claim that the latter can be used instead of the one by Doplicher and Roberts. Of course,
this entails that the integrality assumption on the dimension appears in all statements.
I thank Dr.\ Brugui\`{e}res for correspondence on our -- otherwise independent -- works.

\end{document}